\documentclass[11pt]{article}
\usepackage{latexsym}

\usepackage[top=1in, bottom=1in, left=1in, right=1in]{geometry}

\usepackage{adjustbox}

\usepackage{graphicx} 
\usepackage{color} 
\usepackage{amsmath, amsthm, amssymb}
\usepackage{enumitem}
\setlist[enumerate]{topsep=0pt,itemsep=-1ex,partopsep=1ex,parsep=1ex}

\usepackage{float}
\usepackage{url}
\usepackage{stackrel}
\usepackage{mathrsfs,dsfont}
\usepackage{caption}
\usepackage{subcaption}
\usepackage{wrapfig}
\usepackage{bbm}
\usepackage{bm}
\usepackage{epigraph}
\usepackage{physics}
\usepackage[colorinlistoftodos, shadow]{todonotes} 

\usepackage{mdwlist}
\usepackage{soul}

\usepackage{pifont}
\newcommand{\cmark}{\ding{51}}%
\newcommand{\xmark}{\ding{55}}%

\newtheorem{theorem}{Theorem}[section]
\newtheorem{corollary}[theorem]{Corollary}
\newtheorem{proposition}[theorem]{Proposition}

\newtheorem{example}[theorem]{Example}
\newtheorem{remark}[theorem]{Remark}
\newtheorem{claim}[theorem]{Claim}

\theoremstyle{definition}
\newtheorem{definition}[theorem]{Definition}

\newcommand{\been}{\begin{enumerate}}
\newcommand{\enen}{\end{enumerate}}
\newcommand{\beit}{\begin{itemize}}
\newcommand{\enit}{\end{itemize}}

\def\xrin{\xrightarrow{t \to \infty}}

\def\SS{\mathcal S}
\def\CC{\mathcal C}

\def\GG{\mathcal{G}}

\def\DD{\mathcal{D}}

\def\BB{\mathcal B}

\def\spn{{\rm span}}

\def\la{\leftarrow}

\def\rlas{\rightleftarrows}
\def\ds{\displaystyle}

\def\cbl{\color{blue}}

\newcommand{\R}{\mathbb{R}}

\newcommand{\Z}{\mathbb{Z}}

\usepackage{tikz}    
\usetikzlibrary{shapes,automata,positioning,arrows,fit}
\usepackage{mathtools}  
\usetikzlibrary{snakes}

\usepackage{caption}
\newcommand{\specialcell}[2][c]{\begin{tabular}[#1]{@{}c@{}}#2\end{tabular}}

\numberwithin{equation}{section}

\usepackage{tikz-3dplot}

\usepackage{relsize}

 \tikzset{every node/.style={auto}}
 \tikzset{every state/.style={rectangle, minimum size=0pt, draw=none, font=\normalsize}}
  \tikzset{bend angle=7}
  
  \usepackage{array}
  
  \def\wt{\widetilde}

\usepackage{authblk}

\newcounter{break}

\begin{document}

\title{Foundations of Static and Dynamic\\Absolute Concentration Robustness}

\author[a]{Badal Joshi}
\author[b]{Gheorghe Craciun}

\affil[a]{Department of Mathematics, California State University San Marcos. }
\affil[b]{Departments of Mathematics and Biomolecular Chemistry, University of Wisconsin-Madison.}

\date{}

\maketitle

\begin{abstract}
\noindent Absolute Concentration Robustness (ACR) was introduced by Shinar and Feinberg \cite{shinar2010structural} as robustness of equilibrium species concentration in a mass action dynamical system. 
Their aim was to devise a mathematical condition that will ensure robustness in the function of the biological system being modeled. 
The robustness of function rests on what we refer to as {\em empirical robustness} -- the concentration of a species remains unvarying, when measured in the long run, across arbitrary initial conditions. 
Even simple examples show that the ACR notion introduced in \cite{shinar2010structural} (here referred to as {\em static ACR}) is neither necessary nor sufficient for empirical robustness. 
To make a stronger connection with empirical robustness, we define {\em dynamic ACR}, a property related to long-term, global dynamics, rather than only to equilibrium behavior. 
We discuss general dynamical systems with dynamic ACR properties as well as parametrized families of dynamical systems related to reaction networks. 
We find necessary and sufficient conditions for dynamic ACR in complex balanced reaction networks, a class of networks that is central to the theory of reaction networks.
~ 
\end{abstract}

\section{Introduction}

This work is concerned with the conditions required for {\em \bf empirical robustness} of the concentration of a species in a reaction network. By empirical robustness, we mean that the measured value of a species concentration in the long run remains unchanged even when other conditions, especially the initial concentrations of reagents, change dramatically. 
Shinar, Alon, and Feinberg considered ``{\em the robustness of {\bf equilibrium} species concentrations against fluctuations in the overall reactant supply}'' \cite{shinar2009sensitivity} (bold fonts and italics ours). 
Our goal is to broaden the inquiry by studying the {\bf dynamics of the system}, not {\it merely} the equilibrium values. 
This requires careful consideration of issues related to convergence to a robust value.  
Moreover, since we want to allow arbitrary initial conditions, we must consider global dynamics of the system, a task much more difficult than studying the location of steady states.

We first describe robustness in a biochemical system and the means of its experimental detection via an idealized experimental design. Then we discuss the mathematical property that closely reflects empirical robustness. 

\begin{figure} 
\centering
\includegraphics[scale=0.5]{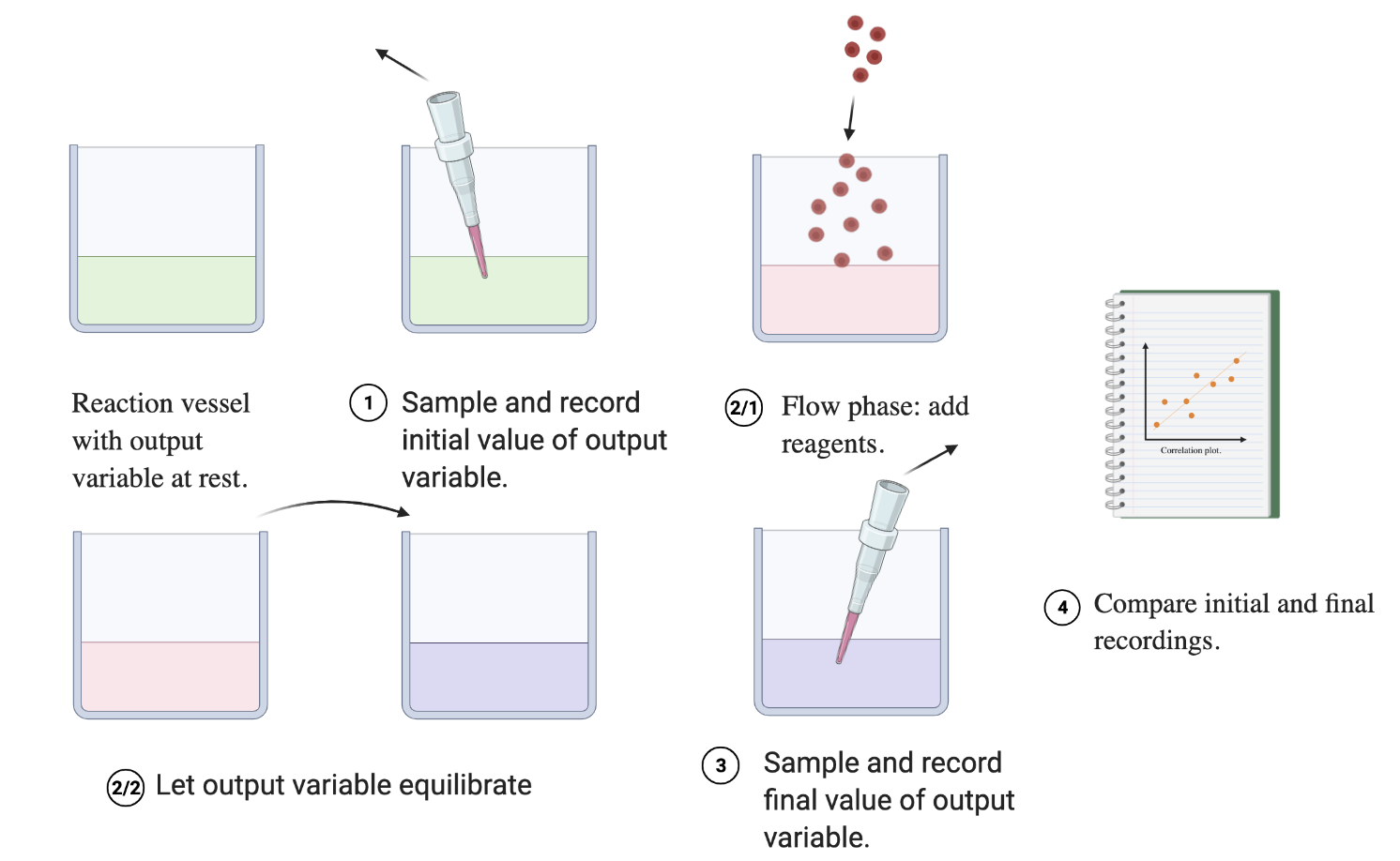}
\caption{A schematic to illustrate the idealized experimental design to identify empirical robustness.}
\label{fig:ryho3uhoehgr}
\end{figure}

\noindent {\bf An idealized experimental design to detect species robustness in a biochemical system} (see Figure \ref{fig:ryho3uhoehgr}).

\begin{figure}[h!] 
\centering
\begin{subfigure}[b]{0.45\textwidth}
\includegraphics[scale=0.41]{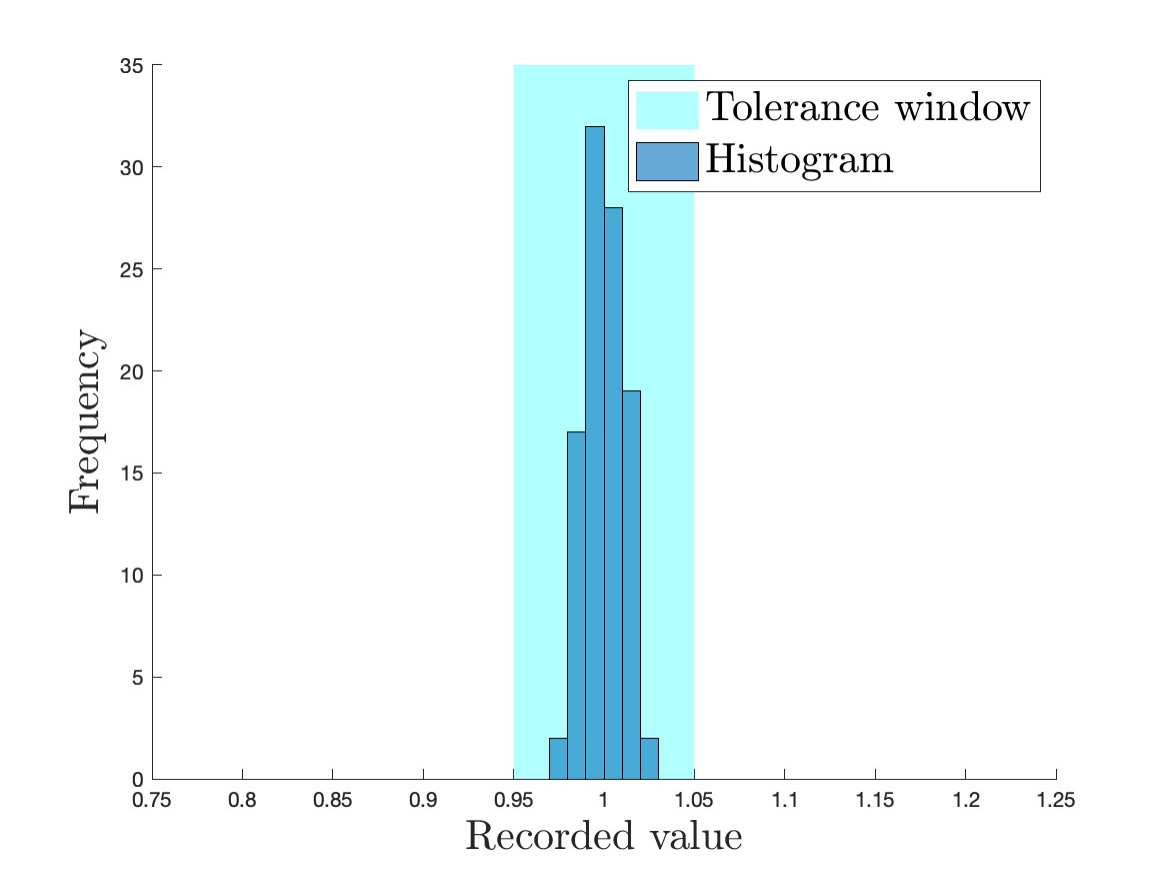}
\caption{Output variable at rest}
\end{subfigure}
\begin{subfigure}[b]{0.45\textwidth}
\includegraphics[scale=0.41]{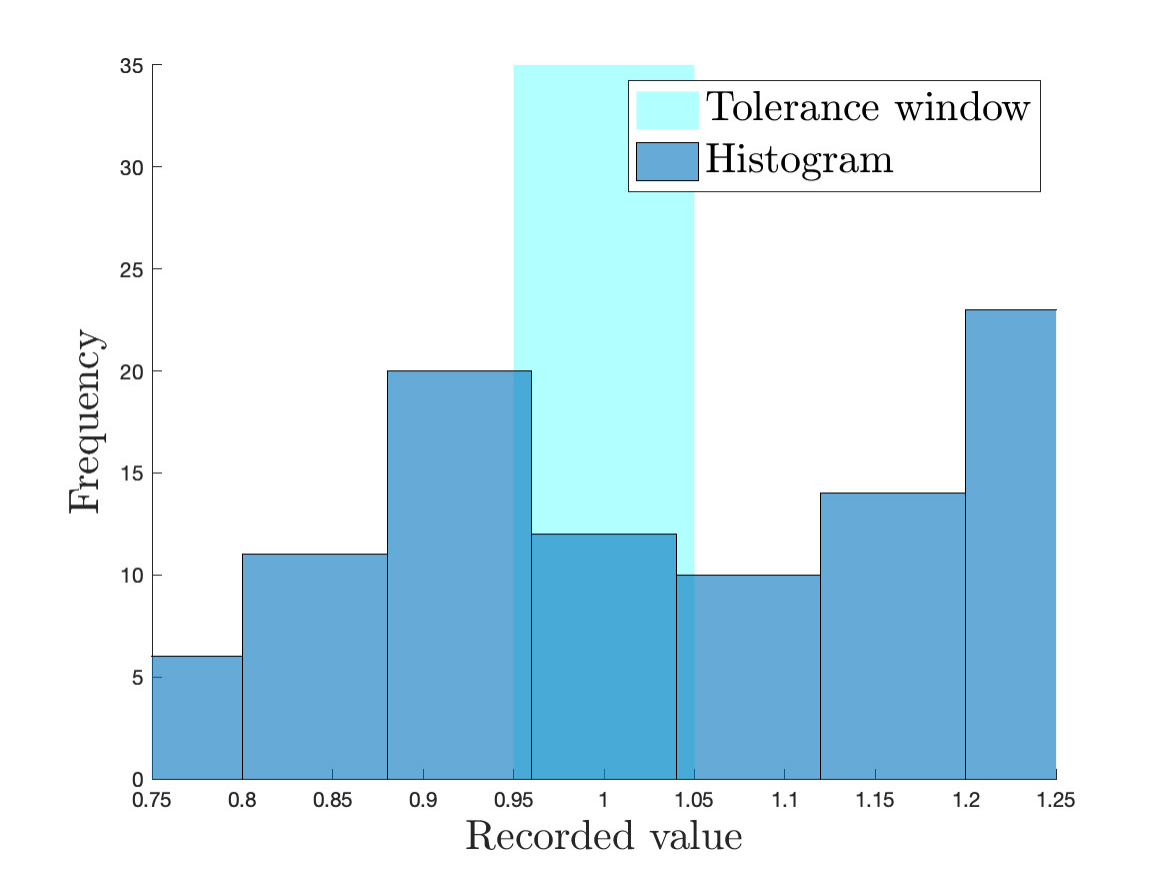}
\caption{Output variable changing with time}
\end{subfigure}
\caption{{\bf Multiple recordings are necessary to distinguish between `output variable at rest' and `output variable changing with time'.} 
Either only initial recordings (before flow) or final recordings (after flow) of the output variable are used. 
Due to small measurement errors, recordings of the output variable at rest may not coincide exactly but should fall within a small tolerance window (left). 
If the output variable is changing in time, we are likely to see a distribution similar to the one on the right.
}
\label{fig:hldfglejhtljeth}
\end{figure}

\begin{figure}[h!] 
\centering
\begin{subfigure}[b]{0.45\textwidth}
\includegraphics[scale=0.41]{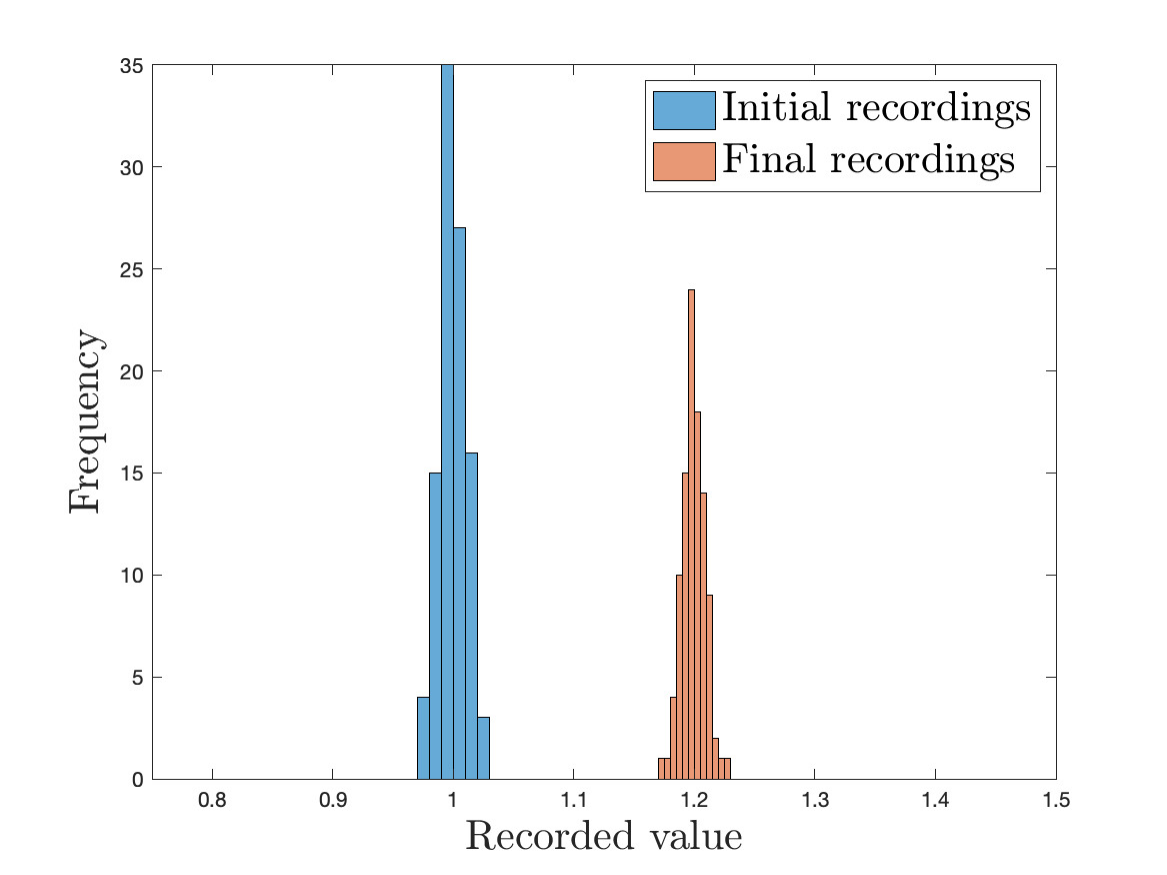}
\caption{Distributions are well-separated: Output variable sensitive to flow process}
\end{subfigure}
\begin{subfigure}[b]{0.45\textwidth}
\includegraphics[scale=0.41]{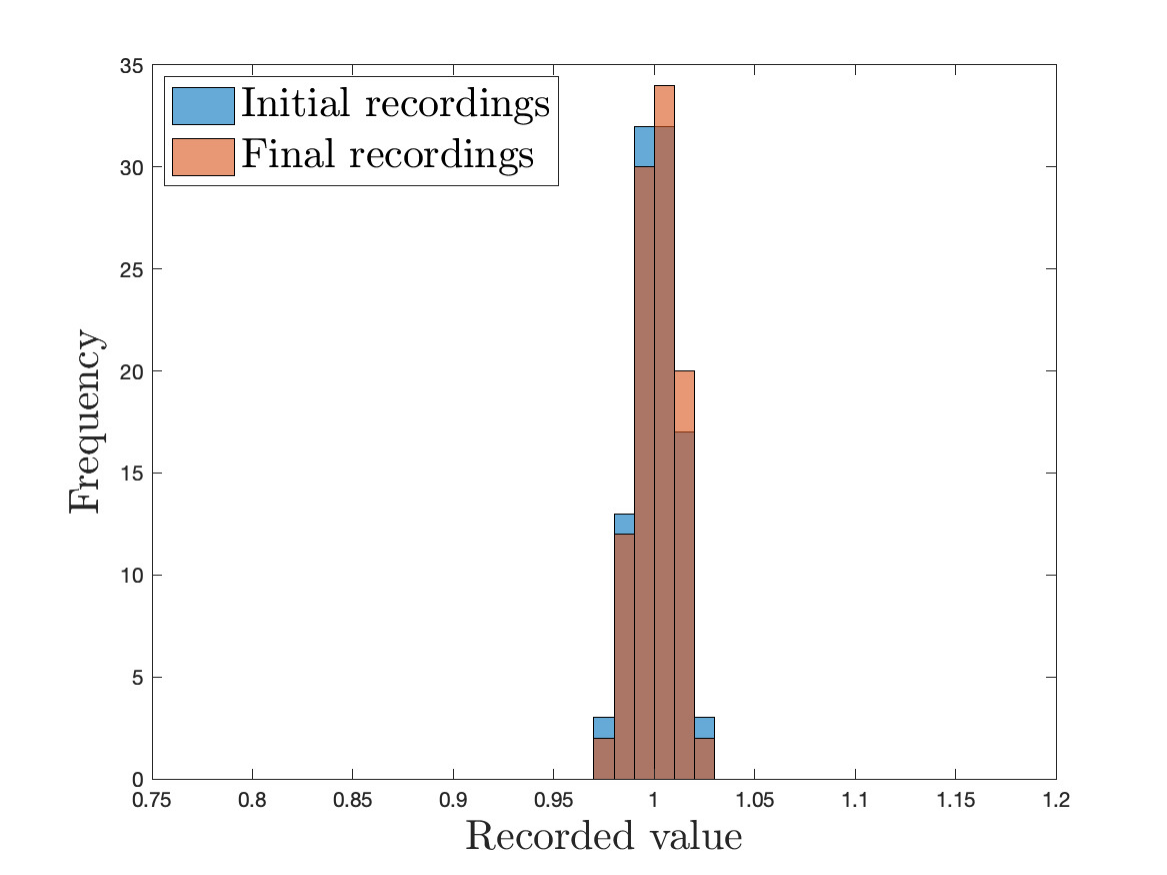}
\caption{Indistinguishable distributions: Output variable has empirical robustness}
\end{subfigure}
\caption{We compare the initial and final recordings of the output variable to find evidence in favor of (right) or against (left) empirical robustness.}
\label{fig:elkrghoewihgoig}
\end{figure}

\been[label={\bf {Step} {\arabic*}.},wide, labelwidth=!, labelindent=0pt]
\item[{\bf {Step} 0.}] {\bf (statistical tolerance):} Designate a variable $\mathfrak{X}$ as an output variable. For concreteness, assume the output variable is the concentration of some biochemical species. Decide an error threshold $p$ for multiple recordings made under identical circumstances. 
\item {\bf (initial recording):} Perform multiple recordings of the output variable $\mathfrak{X}$ at successive time points $t_1, \ldots, t_k$, such that the inter-recording intervals are all different: $t_{i+1} - t_i \ne t_{j+1}-t_j$ for $i \ne j$. 
Denote the set of recordings by $R^I$. 
Sufficiently small variance in the set of initial recordings, $Var[R^I] < p$, provides evidence that $\mathfrak X$ is at rest, and not oscillating or growing or evolving in time (see Figure \ref{fig:hldfglejhtljeth}).  
Multiple recordings also ensure a more accurate estimation of the true rest value of the output variable. 
Denote the average value of the initial recordings by $E[R^I]$. 
\item {\bf (flow phase):} In this phase, a supply of reactants is added to increase the overall concentrations in the reaction volume. The influx can be instantaneous, constant over a large time period, increasing in time, or some other more complicated function of time. 
Eventually, the flow ceases and a large time period is allowed to lapse. 
\item {\bf (final recording):} A second set of recordings of $\mathfrak{X}$ is made. Denote this set of recordings by $R^F$. 
If this set of recordings shows a large variance, $Var[R^F] > p$, then more time is allowed to lapse before this step is repeated. 
If eventually it is the case that $Var[R^F] < p$, then this provides evidence that the output variable has settled down to a rest value. 
Denote the average value of the final recordings by $E[R^F]$.  
\item {\bf (comparison/analysis):} 
The distributions of the initial and the final recordings are compared (see Figure \ref{fig:elkrghoewihgoig}). 
The simplest comparison is that of the average values of the initial and the final recordings.  For instance, suppose that the difference between the two averages is small, i.e. comparable in magnitude to the measurement accuracy and design tolerance: 
\[
\abs{E[R^I] - E[R^F]} \approx \max\left(\sqrt{Var[R^I]}, \sqrt{Var[R^F]}\right).
\]
Further suppose that this finding holds up over repeated trials and a wide variety of flow/influx conditions. Then we conclude that the experiment provides strong evidence in favor of the hypothesis that {\bf the  measured variable shows empirical robustness to the influx process}. 
\enen
The formula in Step 4 is only meant to be suggestive of a statistical test to compare two group means: before-flow group and after-flow group. 
The left side is just the absolute value of the group mean difference. 
The right side is a measure of the within-group variability.
To get some evidence of robustness requires that the between-group difference not be too much higher than the within-group variability.

\noindent{\bf An important remark on the experimental design:} It is worth emphasizing that we are not making any claims about the state of the variables that are not recorded. A non-recorded species concentration may be oscillating, growing in time, converging to zero, or otherwise evolving in time. 

Now we consider issues related to mathematical modeling of a biochemical system. What mathematical property describes empirical robustness sufficiently accurately? 
An important step in the direction of answering this question was taken by Shinar and Feinberg \cite{shinar2010structural}, who defined absolute concentration robustness (ACR). We quote:
\begin{quote}
\tt{A biological system shows absolute concentration robustness (ACR) for an active molecular species if the concentration of that species is identical in every positive steady state the system might admit.}
\end{quote}
Mathematically, the statement is equivalent to: \emph{All positive steady states of the resulting dynamical system are in some hyperplane 
$\{x_i = a_i^*\}$.}  The condition ensures that if $x_i$ is designated as the output variable, then it will remain invariant across positive steady states. 
Shinar and Feinberg \cite{shinar2010structural} then gave a remarkable (sufficient but not necessary) network condition for ACR: `{\em Suppose that a reaction network has deficiency one}, {\em and two non-terminal complexes} (see Section \ref{subsec:backgroundinfo} for definitions) {\em differ in exactly one species. Then the concentration of that species shows ACR.}' The appeal of this criterion is that the network conditions can be checked fairly easily, and they immediately reveal the ACR property for an entire parametrized family of dynamical systems associated with the reaction network. 

Clearly, the elegance and simplicity of the Shinar-Feinberg criterion lends weight to their notion of ACR (which from now we refer to as {\em static ACR}). 
However, if the aim is to model empirical robustness, the definition of static ACR misses the mark to some extent. We go on to quote from \cite{shinar2010structural}: 
\begin{quote}
\tt{The function of an ACR-possessing system is thereby protected even against large changes in the overall supply of the system's components.}
\end{quote}
Clearly, the function of a biochemical system must depend on measurable aspects such as species concentrations. 
We show that {\em static ACR by itself does not confer empirical robustness} and so static ACR is not enough to preserve the function of the system. 
Moreover, there might be systems which lack static ACR, including some that have no steady states whatsoever, but nevertheless the system has a variable that shows empirical robustness. 
We illustrate these points by first giving two simple examples of networks which have static ACR as well as satisfy the Shinar-Feinberg criterion but nevertheless fail to model empirical robustness. Following this, we give two examples of networks which do not have static ACR and yet show robustness in an output variable. 

First, consider the reaction network depicted below in Figure \ref{fig:ergheoruhgo}(a), along with some sample trajectories in Figure \ref{fig:ergheoruhgo}(b) (see \ref{ex8} for more details). 
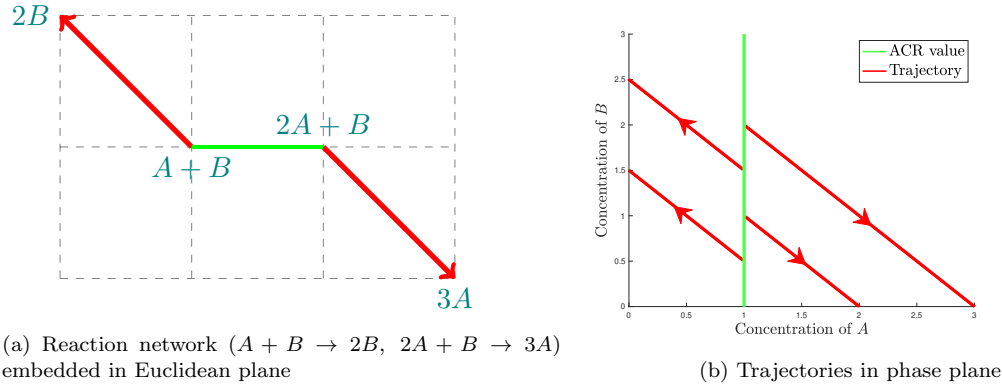
\begin{figure}[h!] 
\centering
\begin{subfigure}[b]{0.45\textwidth}
\begin{tikzpicture}[scale=1.75]
\draw[help lines, dashed, line width=0.25] (0,0) grid (3,2);
\node [below] at (1,1) {{\color{teal} $A+B$}};
\node [left] at (0,2) {{\color{teal} $2B$}};
\node [above] at (2,1) {{\color{teal} $2A+B$}};
\node [below] at (3,0) {{\color{teal} $3A$}};
\draw [->, line width=2, red] (1,1) -- (0,2);
\draw [->, line width=2, red] (2,1) -- (3,0);
\draw [-, line width=1.5, green] (1,1) -- (2,1);
\end{tikzpicture}
\caption{Reaction network ($A+B \to 2B,~ 2A+B \to 3A$) embedded in Euclidean plane}
\end{subfigure}
\begin{subfigure}[b]{0.45\textwidth}
\includegraphics[scale=0.3]{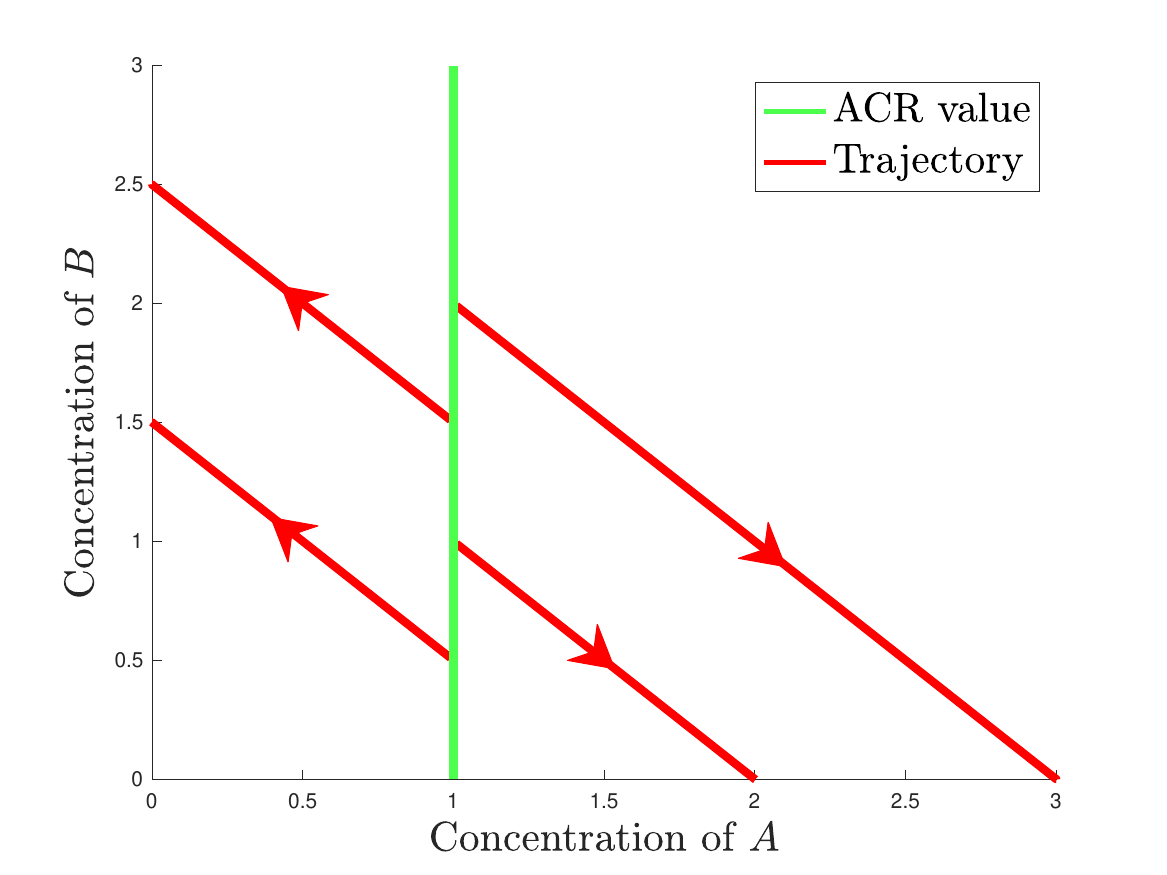}
\caption{Trajectories in phase plane}
\end{subfigure}
\caption{ 
The Euclidean embedding of a reaction network ($A+B \to 2B,~ 2A+B \to 3A$) (left) is important for determining the phase plane trajectories (right). In the one-dimensional case, trajectories (shown in red) are simply parallel to the reaction vectors. The vertical green line in the phase plane is made up of positive steady states. The static ACR property is due to the fact that all positive steady states are contained in a line (hyperplane when there are 2 species) parallel to a coordinate axis. 
However, all positive steady states are unstable, leading to trajectories moving away from ACR value.
Note that this reaction network satisfies the Shinar-Feinberg criterion for static ACR.
}
\label{fig:ergheoruhgo}
\end{figure}
The network in Figure \ref{fig:ergheoruhgo}(a) satisfies the Shinar-Feinberg ACR criterion: 
\beit
\item the deficiency is one -- the two reactions only span a 1 dimensional subspace instead of 2 (see Section \ref{subsec:backgroundinfo}),
\item the non-terminal complexes $A+B$ and $2A+B$ differ by the species $A$ -- the reactant polytope (green line) is parallel to the $A$ axis, 
\enit 
which implies that the concentration of $A$ shows static ACR. Indeed we can see in the figure on the right that all positive steady states lie on a vertical line. 

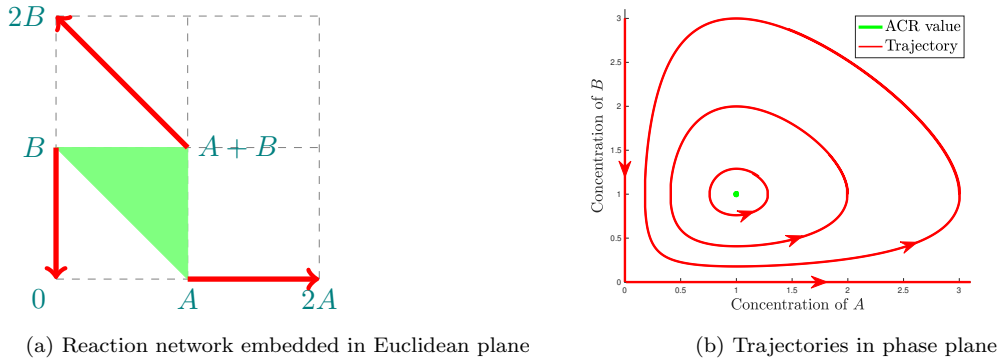
\begin{figure}[h!]
\centering
\begin{subfigure}[b]{0.45\textwidth}
\begin{tikzpicture}[scale=1.75]
\draw[help lines, dashed, line width=0.25] (0,0) grid (2,2);
\node [right] at (1,1) {{\color{teal} $A+B$}};
\node [left] at (0,2) {{\color{teal} $2B$}};
\node [left] at (0,1) {{\color{teal} $B$}};
\node [below left] at (0,0) {{\color{teal} $0$}};
\node [below] at (1,0) {{\color{teal} $A$}};
\node [below] at (2,0) {{\color{teal} $2A$}};

\draw[fill=green!50, draw=green!50, line width=0] (1,1) -- (0,1) -- (1,0) -- cycle;
\draw [->, line width=2, red] (1,1) -- (0,2);
\draw [->, line width=2, red] (0,1) -- (0,0);
\draw [->, line width=2, red] (1,0) -- (2,0);
\end{tikzpicture}
\caption{Reaction network embedded in Euclidean plane}
\end{subfigure}
\begin{subfigure}[b]{0.45\textwidth}
\includegraphics[scale=0.3]{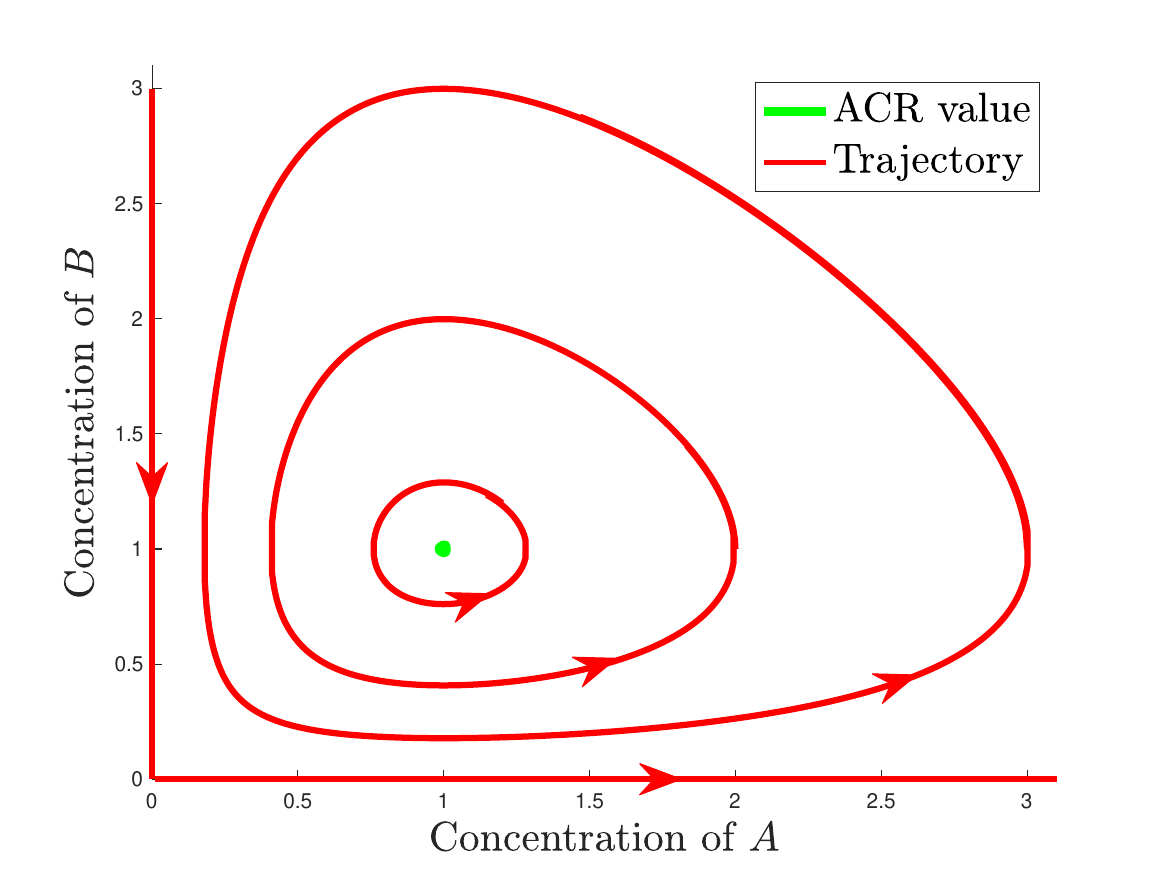}
\caption{Trajectories in phase plane}
\end{subfigure}
\caption{A reaction network ($A \to 2A, A+B \to 2B, B \to 0$) that has static ACR in both species by the Shinar-Feinberg criterion. However, no initial condition leads to convergence to the ACR value.}
\label{fig:erlhkrlhlhllh}
\end{figure}

However, as the numerical solutions in Figure \ref{fig:ergheoruhgo}(b) show, every positive steady state is unstable, and any initial condition (other than the unstable steady state), will result in extinction of one of the species.
The trouble with the last example is that there is another attracting set outside the hyperplane of interest, in this instance a set of boundary steady states. 
Even when there are no other attracting sets, we are not guaranteed convergence to the static ACR hyperplane, as shown by the classic Lotka-Volterra system (Figure \ref{fig:erlhkrlhlhllh}(a), also see \ref{ex10}). 
The network in Figure \ref{fig:erlhkrlhlhllh}(a) satisfies the Shinar-Feinberg ACR criterion: 
\beit
\item the deficiency is one -- the three reactions only span a 2 dimensional subspace instead of 3,
\item the reactant polytope (green triangle) whose vertices are non-terminal complexes has edges parallel to the $A$ axis and the $B$ axis, 
\enit 
which implies that the concentrations of both $A$ and $B$ show static ACR. However, the system has no attractors whatsoever. The unique positive steady state implied by the Shinar-Feinberg criterion is not an attractor, none of the cycles in the figure on the right are attracting (because they are densely/continuously nested), the boundary trajectories are not attracting, and the steady state at the origin is not attracting. Even infinity is not an attractor, since all trajectories with positive initial values remain bounded. This system will generically fail to converge to the ACR value in either coordinate and thus is not a candidate for empirical robustness. 

\begin{figure}
\centering
\begin{subfigure}[b]{0.45\textwidth}
\begin{tikzpicture}[scale=1.75]
\draw[help lines, dashed, line width=0.25] (0,0) grid (1,1);

\node [below left] at (0,0) {{\color{teal} $0$}};
\node [below] at (1,0) {{\color{teal} $A$}};
\node [left] at (0,1) {{\color{teal} $B$}};

\draw [-, line width=1.5, green] (0,0) -- (1,0);
\draw [->, line width=2, red] (0,0) -- (1,0);
\draw [->, line width=2, red] (1,0) -- (0,1);

\end{tikzpicture}
\caption{Reaction network embedded in Euclidean plane}
\end{subfigure}
\begin{subfigure}[b]{0.45\textwidth}
\includegraphics[scale=0.3]{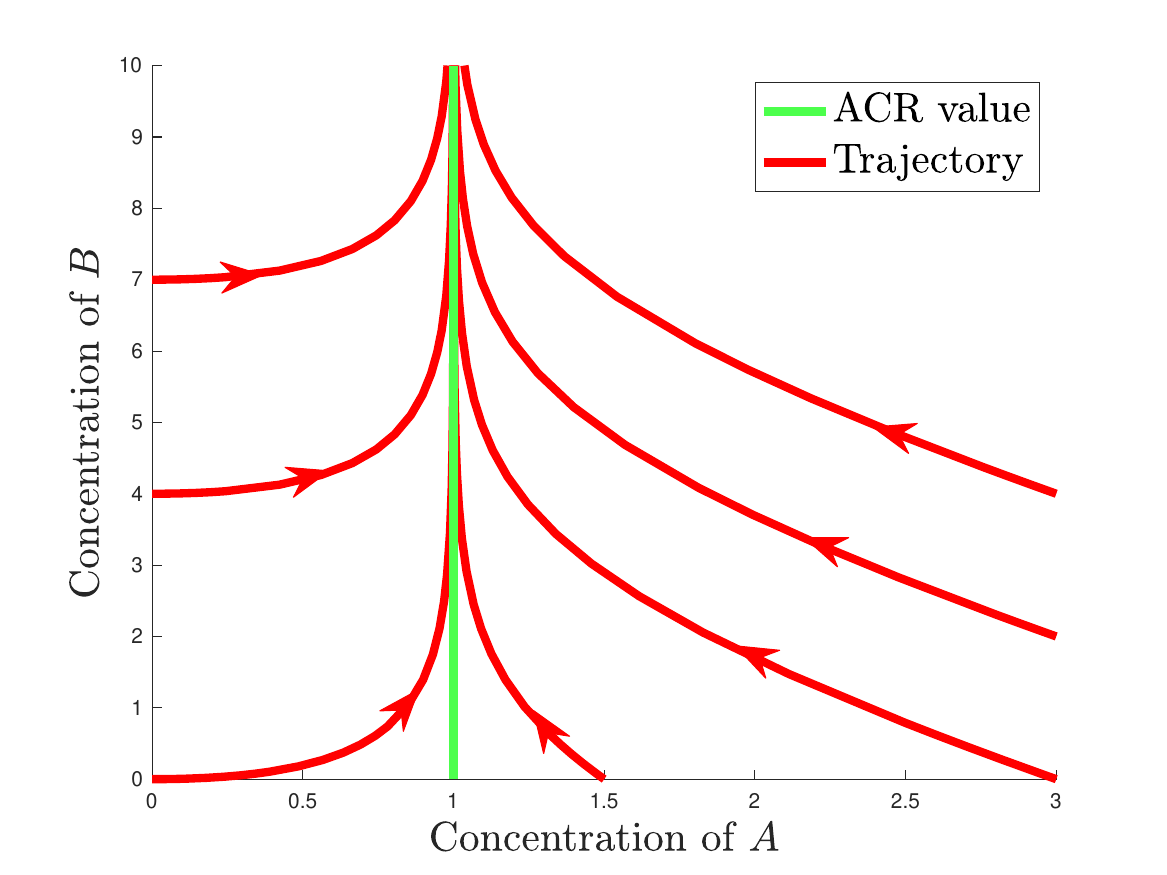}
\caption{Trajectories in phase plane}
\end{subfigure}
\caption{A reaction network ($0 \to A \to B$) that is not ACR and has no steady states, but shows convergence to the same value of concentration of $A$ despite trajectories diverging to infinity.}
\label{fig:4o5uhy4o86hyoii}
\end{figure}
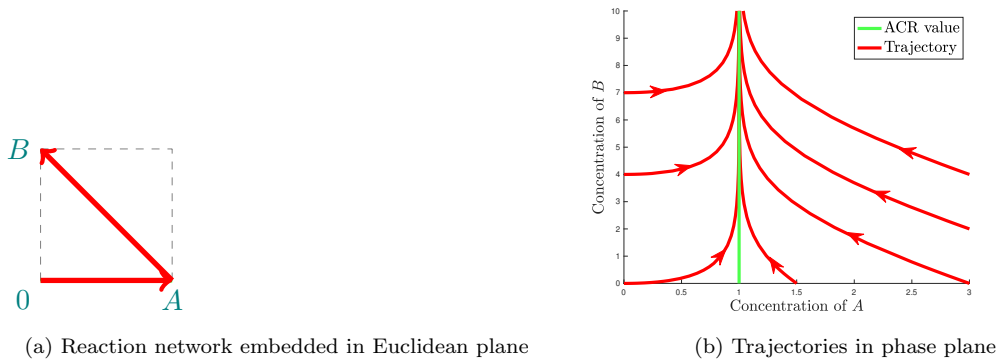

Now we move on to the flip side of the robustness coin. Consider the extremely simple monomolecular network $0 \to A \to B$ shown in Figure \ref{fig:4o5uhy4o86hyoii}(a). $A \to B$ might represent the inactivation of a protein, $A$ being the active form and $B$ being the inactive form, while $0 \to A$ might be a transport process that replenishes the active form. 

This network has deficiency 0, so the Shinar-Feinberg criterion does not apply. Clearly, all trajectories diverge to infinity, so there are no steady states at all. Nevertheless, as shown in Figure \ref{fig:4o5uhy4o86hyoii}(b), all initial conditions result in the concentration of $A$ converging to a unique value. If we designate the concentration of $A$ as the output variable, then the system will show empirical robustness in its value. 

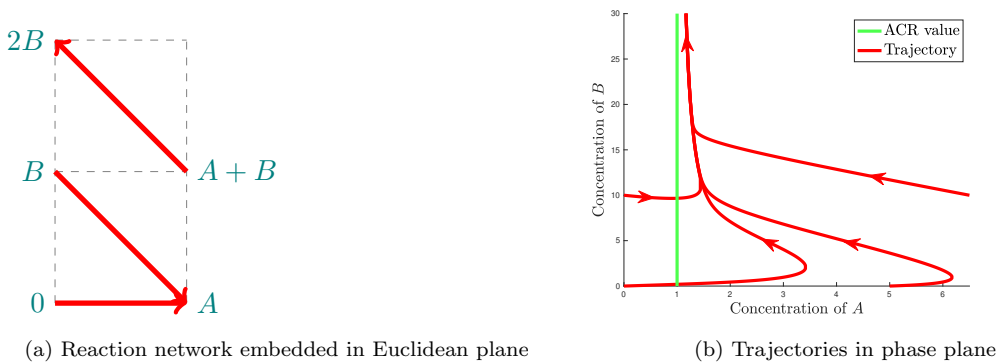
\begin{figure}[h!] 
\centering
\begin{subfigure}[b]{0.45\textwidth}
\begin{tikzpicture}[scale=1.75]
\draw[help lines, dashed, line width=0.25] (0,0) grid (1,2);

\node [left] at (0,0) {{\color{teal} $0$}};
\node [right] at (1,0) {{\color{teal} $A$}};
\node [left] at (0,1) {{\color{teal} $B$}};
\node [right] at (1,1) {{\color{teal} $A+B$}};
\node [left] at (0,2) {{\color{teal} $2B$}};

\draw [->, line width=2, red] (0,1) -- (1,0);
\draw [->, line width=2, red] (1,1) -- (0,2);
\draw [->, line width=2, red] (0,0) -- (1,0);

\end{tikzpicture}
\caption{Reaction network embedded in Euclidean plane}
\end{subfigure}
\begin{subfigure}[b]{0.45\textwidth}
\includegraphics[scale=0.3]{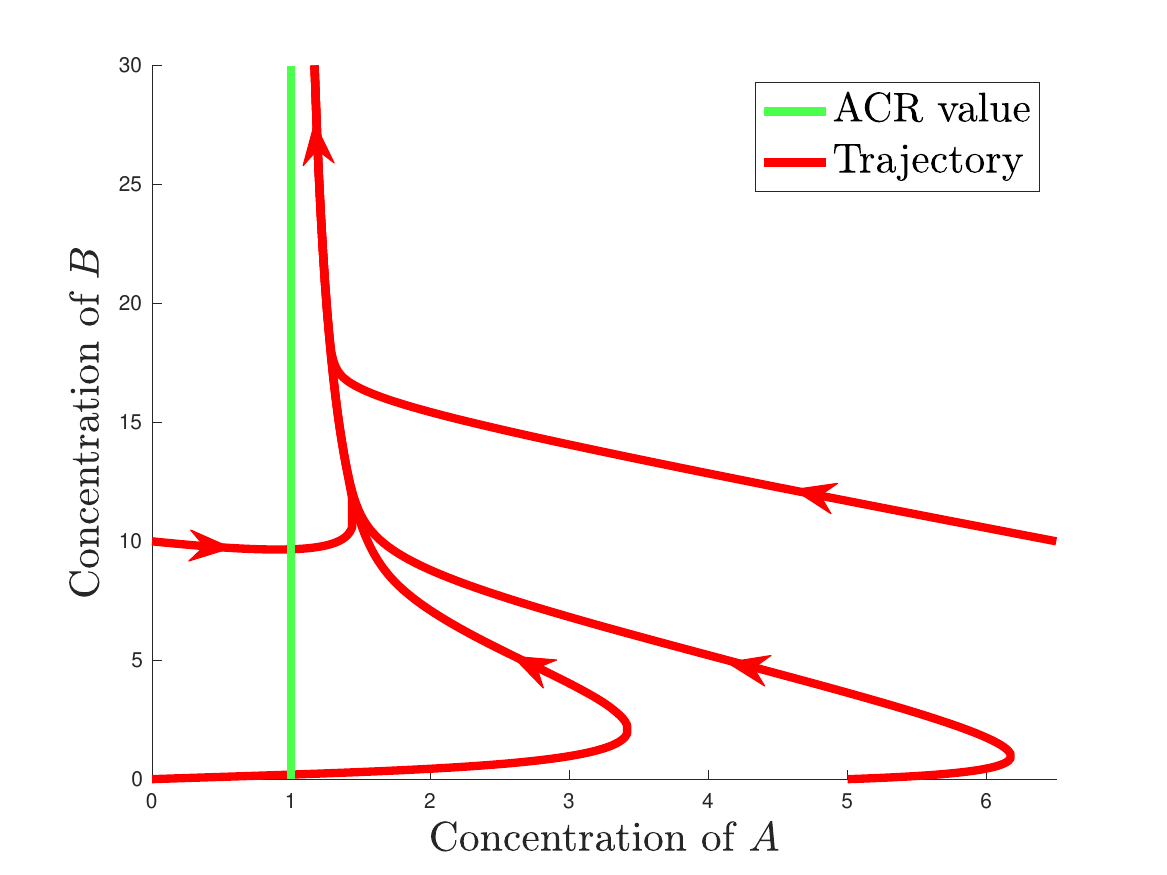}
\caption{Trajectories in phase plane}
\end{subfigure}
\caption{A reaction network ($A+B \to 2B, B \to A, 0 \to A$) that is not ACR and has no steady states, but shows convergence to the same value of concentration of $A$ despite trajectories diverging to infinity.}
\label{fig:45yjhiorthjgejtpiphoi}
\end{figure}

The next example, shown in in Figure \ref{fig:45yjhiorthjgejtpiphoi}(a), is similar to the previous one in that there are no steady states, all trajectories go to infinity, and yet the concentration of one variable converges to a robust value. 
Moreover, unlike the previous example, the robust value does not change as the inflow rate of the robust species $A$ is changed. 
The robust value depends only on the rates of the true chemical reactions $A+B \to 2B$ and $B \to A$, and not on the rate of the transport/inflow reaction $0 \to A$. A rigorous analysis of this system will appear in future work, here we show the robustness in concentration of $A$ by simulating some trajectories, shown in Figure \ref{fig:45yjhiorthjgejtpiphoi}(b). 

We now return to the project of constructing a theoretical framework that would better capture the property of empirical robustness in the context of mathematical models based on deterministic dynamical systems.
{\bf The simplest way to proceed seems to be to insist that the hyperplane $\{x \in \R^n_{\ge 0} ~|~ x_i = a_i^* > 0\}$ as a whole be an attractor to all initial conditions that are compatible with the hyperplane} (see Definition \ref{def:6ijeqjgqjeh'ije}).
We will call this notion {\em dynamic ACR}.
Clearly dynamic ACR requires that there are no attracting sets outside $\{x_i = a_i^*\}$. But it does not require that $\{x_i = a_i^*\}$ be invariant, see for instance the network in Figure \ref{fig:45yjhiorthjgejtpiphoi}(a)--the figure shows trajectories (red curves) crossing over the attracting hyperplane (green line) but eventually converging to it.

Establishing clear mathematical foundations for the study of empirical robustness is essential for the theory to make consequential predictions relevant to biochemistry. 
Empirical robustness has been observed experimentally in a large class of bacterial two-component signaling systems \cite{russo1993essential,hsing1998mutations,batchelor2003robustness,shinar2007input}. 
The circuit design for signal transduction, where a signal is transported from the cell environment to its interior, uses a mechanism involving a bifunctional component \cite{alon2019introduction}. 
A bifunctional component exerts two opposing forces, for instance promoting phosphorylation as well as dephosphorylation of a substrate. 
Such a mechanism ensures that the output depends on the signal strength but not on the details of the circuit implementation, for instance the number of signaling proteins that form the circuit. 
In future work, we will prove that networks with a bifunctional component have the property of dynamic ACR, and not merely static ACR. 
In particular this means that in a signal transduction circuit with a bifunctional component, for any initial value (which encodes the circuit implementation), the cell response converges to a value that only depends on the signal strength.

This article is organized as follows.  
Section~\ref{sec:q3o8thegqpw} contains the central definitions of this article, and some propositions to highlight the connections with existing notions. 
Section~\ref{subsec:backgroundinfo} contains some background information on deterministic modeling of the dynamics of reaction networks and previous work on static ACR.
Section~\ref{sec:illustrativeexamples}, the main course, has several illustrative examples which delineate the specific conditions in the definitions.
Section~\ref{sec:3pyihoeqhljd} applies the previous concepts to reaction networks and contains a discussion on static and dynamic ACR at the network level.
Section~\ref{sec:5yiu6fikfyuerlkgfye} is about the interplay between ACR and complex balance, two important ideas in reaction network theory.

\section{Basic Definitions of Static ACR and Dynamic ACR in real dynamical systems} \label{sec:q3o8thegqpw}

Throughout the article, we consider a dynamical system $\DD$ defined by $\dot x = f(x)$ with $x \in \R^n_{\ge 0}$ and a  smooth vector field $f$ for which $\R^n_{\ge 0}$ is forward invariant. A point $x_0 \in \R^n_{\ge 0}$ is a {\em steady state} of $\DD$ if $f(x_0)=0$. 

\begin{definition} \label{def:0496yuwpeh}
The {\em kinetic subspace} of $\DD$ is defined to be the linear span of the image of $f$, denoted by $\spn(\Im(f))$. 
The points $x,y \in \R^n_{\ge 0}$ are {\em compatible} if $y - x \in \spn(\Im(f))$. 
The sets $S, S' \subseteq \R^n_{\ge 0}$ are {\em compatible} if there are $x \in S$ and $x' \in S'$ such that $x$ and $x'$ are compatible.
A {\em compatibility class $S$} is a nonempty subset of $\R^n_{\ge 0}$ such that $x,y \in S$ if and only if $y - x \in \spn(\Im(f))$.
\end{definition}
We first define the Shinar-Feinberg notion of ACR, which we refer to as static ACR since it relates only to a property of the set of steady states and since it allows us to make a distinction with dynamic ACR.
\begin{definition} \label{def:42u5hgd;;goi}
$\DD$ is a {\em static ACR system} if $\DD$ has a positive steady state and there is an $i \in \{1,\ldots, n\}$ and a positive $a_i^* \in \R_{> 0}$ such that any positive steady state $x \in \R^n_{> 0}$ is contained in the hyperplane $\{x_i = a_i^*\}$. 
Any such $x_i$ and $a_i^*$ is a {\em static ACR variable} and its {\em static ACR value}, respectively.
\end{definition}
We now introduce dynamic ACR with the explicit goal of more accurately modeling empirical robustness. 

\begin{definition} \label{def:6ijeqjgqjeh'ije}
$\DD$ is a {\em dynamic ACR system} if there is an $i \in \{1,\ldots, n\}$ with $f_i \not \equiv 0$ and a positive $a_i^* \in \R_{> 0}$ such that for any $x(0) \in \R^n_{> 0}$ that is compatible with $\{x \in \R^n_{>0} ~|~ x_i = a_i^*\}$, a unique solution to $\dot x = f(x)$ exists up to some maximal $T_0(x(0)) \in (0, \infty]$, and $x_i(t) \xrightarrow{t \to T_0} a_i^*$. Any such $x_i$ and $a_i^*$  is a {\em dynamic ACR variable} and its {\em dynamic ACR value}, respectively.  
\end{definition}
If the dynamical system $\dot x = f(x)$ does not have the possibility of a finite-time blow-up, then $T_0(x(0)) = \infty$ for any $x(0) \in \R^n_{>0}$. 
None of the examples in this paper have the possibility of a finite-time blow-up, however the definition is more generally applicable to allow for this possibility. The use of ``its'' dynamic ACR value in Definition \ref{def:6ijeqjgqjeh'ije} is justified by the next result. 
\begin{theorem} \label{thm:oeh13othop3h}
Suppose that $\DD$ is a static (dynamic resp.) ACR system. Each static (dynamic resp.) ACR variable has a unique static (dynamic resp.) ACR value.  
\end{theorem}
\begin{proof}
The statement about a static ACR system follows immediately from the definition. Suppose that $x_i$ is a dynamic ACR variable with distinct ACR values $a_i^*$ and $b_i^*$. 
Then $\{x \in \R^n_{>0} ~|~ x_i = a_i^*\}$ and $\{x \in \R^n_{>0} ~|~ x_i = b_i^*\}$ are not compatible, which implies that $\{x \in \R^n_{>0} ~|~ x_i = a_i^*\}$ and $\{x \in \R^n_{>0} ~|~ x_i = c_i\}$ are not compatible for any positive $c_i \ne a_i^*$.
But then the set of points compatible with $\{x \in \R^n_{>0} ~|~ x_i = a_i^*\}$ is contained in $\{x \in \R^n_{>0} ~|~ x_i = a_i^*\}$ which implies that $f_i \equiv 0$, a contradiction. So a dynamic ACR variable must have a unique dynamic ACR value. 
\end{proof}

Dynamic ACR requires all compatible initial values to result in convergence of the ACR variable to the ACR value. 
But the set of compatible initial values can be quite different for different systems.

\begin{definition} \label{def:4;otiyj;oi4hg;o}
Suppose that $\DD$ is a dynamic ACR system which has a dynamic ACR variable $x_i$ with ACR value $a_i^*$. 
Let $\pi_i(y)$ denote the projection of $y \in \R^n$ on the $i$-axis and let $N_{i, a_i^*} \coloneqq \{ \pi_i(y) ~|~ y\in \R^n_{> 0} \mbox{ not compatible with } \{x \in \R^n_{>0} ~|~ x_i = a_i^*\} \}.$
The variable $x_i$ is a {\em wide basin dynamic ACR variable} if 
$N_{i, a_i^*}$ has an upper bound. 
Otherwise, $x_i$ is a {\em narrow basin dynamic ACR variable}. 
Finally, $x_i$ is a {\em full basin dynamic ACR variable} if 
$N_{i, a_i^*} = \varnothing$. 
\end{definition}

\begin{theorem} \label{thm:q035yhtihgje}
The following hold for a dynamical system $\DD: \dot x = f(x)$, $x \in \R^n_{\ge 0}$. 
\been
\item $\DD$ has a unique positive steady state if and only if every variable is a static ACR variable.  
\item If $\DD$ has a globally attracting positive steady state (i.e. the basin of attraction is the positive orthant) then every variable is a full basin dynamic ACR variable. 
\item If every variable is a dynamic ACR variable then $\DD$ has a globally attracting positive steady state. 
\enen
\end{theorem}
\begin{proof}
The first two statements are immediate from the definitions. 

In order to prove the third statement, assume that for every $i \in \{1,\ldots,n\}$, $x_i$ is a dynamic ACR variable with dynamic ACR value $a_i^*$. 
Note that, since $f_i \not \equiv 0$, it follows that there exists a neighborhood $N^*$ of $a^* \coloneqq (a_1^*, \ldots, a_n^*)$ in $\R^n_{> 0}$ that is compatible with the hyperplane $\{x_i = a_i^*\}$ for all $i$. Indeed, $f_i \not \equiv 0$ implies that the kinetic subspace of $\DD$ is transversal to $\{x_i = a_i^*\}$, which implies that there exists a neighborhood $N_i^*$ of $a^*$ in $\R^n_{> 0}$ that is compatible with the hyperplane $\{x_i = a_i^*\}$; we can then define 
$$N^* = \bigcap_{1\le i \le n}N_i^*.$$ 

Then, since $x_i$ is dynamic ACR for all $i$, we conclude that for any $x_0 \in N^*$ the trajectory that starts at $x_0$ converges to $a^*$. This allows us to prove that the kinetic subspace of $\DD$ is the whole $\R^n$. Indeed, assume that this is {\em not} true, in order to obtain a contradiction. It follows that there exists some compatibility class (i.e., shifted version of the kinetic subspace, of the form $x_0 + span(Im(f))$) that intersects $N^*$ but does not contain $a^*$; note also that $(x_0 + span(Im(f))) \cap \R^n_{\ge 0}$ is an invariant set of $\DD$. Then it follows that there exist a point $x^* \in N^*$ such that the trajectory that starts at $x^*$ does not converge to $a^*$, a contradiction.

Therefore, the kinetic subspace of $\DD$ is  $\R^n$, which implies that $a^*$ is globally attracting.
\end{proof}

\begin{corollary}
If every variable in a dynamical system $\DD$ is dynamic ACR then every variable in $\DD$ is full basin dynamic ACR.
\end{corollary}

Under some mild additional hypotheses (existence of steady states and compatibility conditions), dynamic ACR implies static ACR for a given variable. 
\begin{theorem}
Consider a dynamical system $\DD$ where $x_i$ is a dynamic ACR variable with ACR value $a_i^*$.
Let $\BB$ denote the set of positive steady states of $\DD$. 
The following are equivalent:
\been
\item $x_i$ is a static ACR variable with static ACR value $a_i^*$. 
\item $\varnothing \ne \BB \subseteq \{y\in \R^n_{> 0}: y \mbox{ compatible with } \{x_i = a_i^*\} \}$. 
\enen
\end{theorem}
\begin{proof}
We first show that ({\it 1} $\implies$ {\it 2}). Suppose that $\BB = \varnothing$. Then $\DD$ is not static ACR and there are no static ACR variables. 
If there is a positive steady state that is not compatible with $\{x \in \R^n_{>0} ~|~ x_i = a_i^*\}$ then in particular, there is a positive steady state which is not on the hyperplane $\{x \in \R^n_{>0} ~|~ x_i = a_i^*\}$, which shows that $\DD$ is not static ACR.  

Now we show that ({\it 2} $\implies$ {\it 1}). 
Suppose $\DD$ has positive steady states and each of these is compatible with $\{x \in \R^n_{>0} ~|~ x_i = a_i^*\}$. 
Consider one such positive steady state, say $z$. By definition, if $(y(t))_{y \ge 0}$ is a trajectory with $y(0) = z$, then $y(t) =z$ for all $t \ge 0$. But since, by definition of dynamic ACR, $y(t) \xrin \{x_i = a_i^*\}$, we must have $z \in \{x_i = a_i^*\}$, i.e. $z_i = a_i^*$. Therefore, $x_i$ is the static ACR variable with static ACR value $a_i^*$. 
\end{proof}

\section{Background information on reaction networks} \label{subsec:backgroundinfo}
The definitions and claims appearing thus far have been about general real dynamical systems. We mostly work with reaction networks and mass action systems, for which we use standard notation and terminology. Here we only give a quick summary of the conventions, see for instance \cite{joshi2015survey} for further details. In Example \ref{ex:reactionnetexample}, we illustrate all concepts defined below. 

Throughout this paper, we use upper case letters ($X,Y, Z, A, B$) for species participating in reactions and the corresponding lower case letters ($x,y,z,a,b$) for their concentrations, which are dynamic, time-dependent quantities. An example of a reaction is $X + Y \to 2Z$, where $X+Y$ is referred to as the source complex, while $2Z$ is the product complex. The rate of any given reaction is a nonnegative-valued function of species concentrations. We usually use mass action kinetics wherein the rate is proportional to the product, taken with multiplicity, of reactant concentrations. The proportionality constant, called the reaction rate constant, is placed adjacent to the reaction arrow, as follows: $X + Y \xrightarrow{k} 2Z$. The rate of this reaction under mass action kinetics is $kxy$. 
The reaction vector for this reaction is the difference between the product complex and the source complex, i.e. $2Z-(X+Y)$, which under a choice of standard basis can also be written as $(-1,-1,2)$. 
A {\em reaction network} is a nonempty set of reactions, such that every species participates in at least one reaction, and none of the reaction vectors is the zero vector. The {\em stoichiometric subspace} of a reaction network is the subspace spanned by the set of reaction vectors of the reaction network.
A reaction network $\GG$ is said to be {\em mass conserving} if there is a positive, linear conservation law involving all species, in other words, if there is a positive vector orthogonal to the stoichiometric subspace. 

We say that two complexes are in the same {\em linkage class} if there is a sequence of reactions (backward or forward) connecting the two complexes. 
For kinetic systems of reaction networks where each linkage class has precisely one terminal strong linkage class  (see Definitions 8, 9, 10 and Theorem in Section 6 of \cite{feinberg1977chemical}), the kinetic subspace in Definition \ref{def:0496yuwpeh} coincides with the {\em stoichiometric subspace}. 
The {\em deficiency} of a reaction network is $\delta = n - \ell - s$, where $n$ is the number of complexes in the reaction network, $\ell$ is the number of linkage classes and $s$ is the dimension of the stoichiometric subspace. 

We use $\GG$ to denote a reaction network and $K$ to denote a specific choice of mass action kinetics for $\GG$, so that $(\GG,K)$ is a mass action dynamical system. 
A mass action system $(\GG,K)$ is complex balanced if at every positive steady state, for each complex $\CC$, the sum of reaction rates where $\CC$ is the reactant complex equals the sum of reaction rates where $\CC$ is the product complex. 
A network is {\em weakly reversible} if every reaction is part of a cycle of reactions. The mass action system $(\GG,K)$ is complex balanced for any choice of $K$ if $\GG$ is weakly reversible and has zero deficiency. 

A complex is {\em non-terminal} if it is not in a terminal strong linkage class. 
\begin{theorem}[Shinar \& Feinberg \cite{shinar2010structural} criterion for static ACR] \label{thm:shinar-feinberg}
Consider a reaction network $\GG$ such that (i) the deficiency of $\GG$ is $1$, and (ii) there are two non-terminal complexes $\CC_1$ and $\CC_2$ in $\GG$ such that $\CC_1 - \CC_2 = \alpha X$ for some $\alpha \ne 0$. Then for any choice of $K$ such that $(\GG,K)$ has a positive steady state, the concentration of $X$ is a static ACR variable in $(\GG,K)$. 
\end{theorem}

\begin{example} \label{ex:reactionnetexample}
An example of a reaction network is 
\begin{align*}
S + E &\stackrel[k_2]{k_1}{\rightleftarrows} C \xrightarrow{k_3} P + E \\
P &\xrightarrow{k_4} S
\end{align*}
The species $\{S,E,C,P\}$ have time-dependent concentrations $\{s(t), e(t), c(t), p(t)\}$, respectively. The reaction $S+E \xrightarrow{k_1} C$ has source complex $S+E$, product complex $C$, mass action reaction rate constant $k_1$ and the mass action reaction rate $k_1se$. 
Assuming an arbitrary ordering of the species set $(S,E,C,P)$, the stoichiometric subspace is a subspace of $\R^4$ spanned by the following set of four reaction vectors (ordered according to their reaction rate constants)
\[
\left \{\begin{pmatrix} -1 \\ -1 \\ 1 \\ 0 \end{pmatrix}, \begin{pmatrix} 1 \\ 1 \\ -1 \\ 0 \end{pmatrix}, \begin{pmatrix} 0 \\ 1 \\ -1 \\ 1 \end{pmatrix}, \begin{pmatrix} 1 \\ 0 \\ 0 \\ -1 \end{pmatrix} \right \}. 
\]
The reaction network is mass conserving since $(1,1,2,1)$ is a positive vector that is orthogonal to the stoichiometric subspace.  In the linkage class $P \to S$, the terminal strong linkage class is $\{S\}$. The reaction network is not weakly reversible since the terminal strong linkage class does not coincide with the linkage class. The reaction network has two linkage classes and each linkage class has precisely one terminal strong linkage class. This implies that the kinetic subspace is same as the stoichiometric subspace. 
The deficiency is $\delta = n - \ell - s = 5 - 2 - 2 = 1$. 
\end{example}

\noindent{\bf Previous work on ACR:}
Before proceeding with the remainder of the paper, we give a brief, and by no means exhaustive, survey of existing literature on static ACR. 
Since the seminal work by Shinar and Feinberg in 2010 \cite{shinar2010structural}, ACR has generated tremendous interest and enthusiasm. 
Shinar and Feinberg gave further results on connections between network structure and ACR properties \cite{shinar2011design}. 
Karp, P{\' e}rez Mill{\' a}n, Dasgupta, Dickenstein, Gunawardena studied the ACR conditions from a broader point of view of complex-linear invariants \cite{karp2012complex}. 
Dexter and Gunawardena \cite{dexter2013dimerization} showed that ``{\em homodimerization of IDH and bifunctionality of its regulatory enzyme}'' lead to robustness in a biochemically realistic mathematical model of the IDH system. 
Dexter, Dasgupta and Gunawardena \cite{dexter2015invariants} gave other classes of invariants besides ACR to include bounds on concentration, hybrid robustness, and robust concentration ratio. 
Stochastic (continuous-time Markov chain) models of reaction networks with the ACR property were studied in \cite{anderson2014stochastic,anderson2017finite,enciso2016transient} and control theory aspects in \cite{cappelletti2020hidden,kim2020absolutely}.
Pascual-Escudero and Feliu \cite{pascual2020local} make a distinction between networks with ACR and a broader class with the property of zero sensitivity with respect to initial conditions.

\section{Illustrative Examples of Static and Dynamic ACR systems} \label{sec:illustrativeexamples}

\subsection{Examples with one species.}

\begin{enumerate}[label={\bf Ex \arabic*}.,wide, labelwidth=!, labelindent=0pt]

\item\label{ex3} {\em (Unique positive steady state implies static ACR.)} 
Consider the network $\{2A \xrightarrow{k'} 3A, ~A \xrightarrow{k} 0\}$ whose mass action system is $\dot a = a(k'a-k)$. Existence of the unique positive steady state $a=k/k'$ implies  the system is static ACR for positive $k$ and $k'$; $a$ is the static ACR variable with value $k/k'$. The steady state is repelling for any choice of $k$ and $k'$, so the system is not dynamic ACR.

\item\label{ex4} {\em (Globally attracting steady state implies dynamic ACR.)} 
Consider the network $0 \stackrel[k']{k}{\rlas} A$ and its mass action system $\dot a = k - k'a$, which has a globally attracting positive steady state $a = k/k'$. Therefore, the system is both static and dynamic ACR. Moreover, $a$ is a full basin dynamic ACR variable with its (static and dynamic) ACR value $k/k'$. 
\setcounter{break}{\value{enumi}}
\end{enumerate}

\subsection{Two species, one-dimensional system, infinitely many steady states.}

\begin{figure}[h!] 
\centering
\begin{subfigure}[b]{0.3\textwidth}
\begin{tikzpicture}[scale=1]
\draw[help lines, dashed, line width=0.25] (0,0) grid (3,2);
\draw [<->, line width=1.5, blue] (0,2.5) -- (0,0) -- (3.5,0);

\node [below] at (2,0) {{\cbl $A$}};
\node [left] at (0,2) {{\cbl $B$}};

\draw [->, line width=2, red] (1,1) -- (0,2);
\draw [->, line width=2, red] (2,1) -- (3,0);

\draw [-, line width=2, green] (1,1) -- (2,1);

\end{tikzpicture}
\caption{\bf Static (but not dynamic) ACR network}
\end{subfigure}
\begin{subfigure}[b]{0.3\textwidth}
\begin{tikzpicture}[scale=1]
\draw[help lines, dashed, line width=0.25] (0,0) grid (2,2);
\draw [<->, line width=1.5, blue] (0,2.5) -- (0,0) -- (2.5,0);

\node [below] at (2,0) {{\cbl $A$}};
\node [left] at (0,2) {{\cbl $B$}};

\draw [->, line width=2, red] (0,1) -- (0,2);
\draw [->, line width=2, red] (1,1) -- (1,0);

\draw [-, line width=2, green] (0,1) -- (1,1);

\end{tikzpicture}
\caption{\bf Static (but not dynamic) ACR network}
\end{subfigure}
\begin{subfigure}[b]{0.3\textwidth}
\begin{tikzpicture}[scale=1]
\draw[help lines, dashed, line width=0.25] (0,0) grid (2,2);
\draw [<->, line width=1.5, blue] (0,2.5) -- (0,0) -- (2.5,0);

\node [below] at (2,0) {{\cbl $A$}};
\node [left] at (0,2) {{\cbl $B$}};

\draw [->, line width=2, red] (1,1) -- (0,2);
\draw [->, line width=2, red] (0,1) -- (1,0);

\draw [-, line width=2, green] (0,1) -- (1,1);

\end{tikzpicture}
\caption{\bf Dynamic ACR network}
\end{subfigure}

\caption{The red arrows depict reactions. An arrow originates at the source complex and terminates at the product complex. {\bf (Left:)} The reaction network $A+B \to 2B, 2A + B \to 3A$. {\bf (Middle:)} The reaction network $A+B \to A, B \to 2B$.  {\bf (Right:)} The reaction network $A+B \to 2B, B \to A$.The green line-segment that joins the two sources complexes depicts the reactant polytope of the network. Since the reactant polytope is parallel to a coordinate axis (the $A$ axis) in all cases, all networks are static ACR, with $A$ as the only static ACR species. Only the network in (c) is dynamic ACR.}
\label{fig:34p98hwrtuhw}
\end{figure}
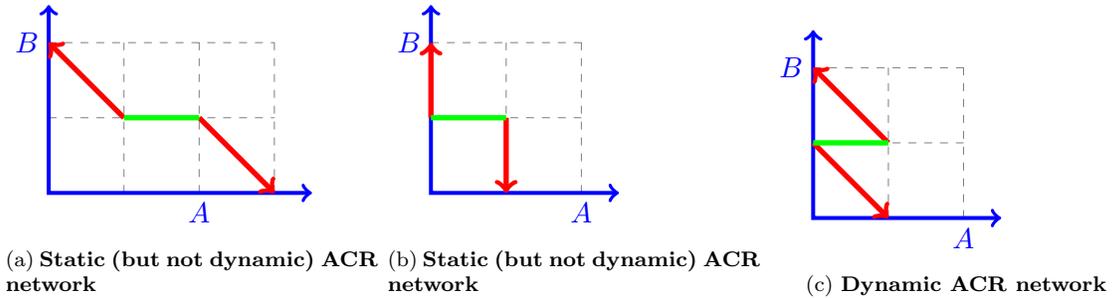

\begin{figure}[h!] 
\centering
\begin{subfigure}[b]{0.45\textwidth}
\begin{tikzpicture}[scale=1.75]
\draw[help lines, dashed, line width=0.25] (0,0) grid (1,2);

\node [right] at (1,0) {{\color{teal} $A$}};
\node [left] at (0,1) {{\color{teal} $B$}};
\node [right] at (1,1) {{\color{teal} $A+B$}};
\node [left] at (0,2) {{\color{teal} $2B$}};

\draw [-, line width=1.5, green] (0,1) -- (1,1);
\draw [->, line width=2, red] (0,1) -- (1,0);
\draw [->, line width=2, red] (1,1) -- (0,2);

\end{tikzpicture}
\caption{Reaction network embedded in Euclidean plane}
\end{subfigure}
\begin{subfigure}[b]{0.45\textwidth}
\includegraphics[scale=0.3]{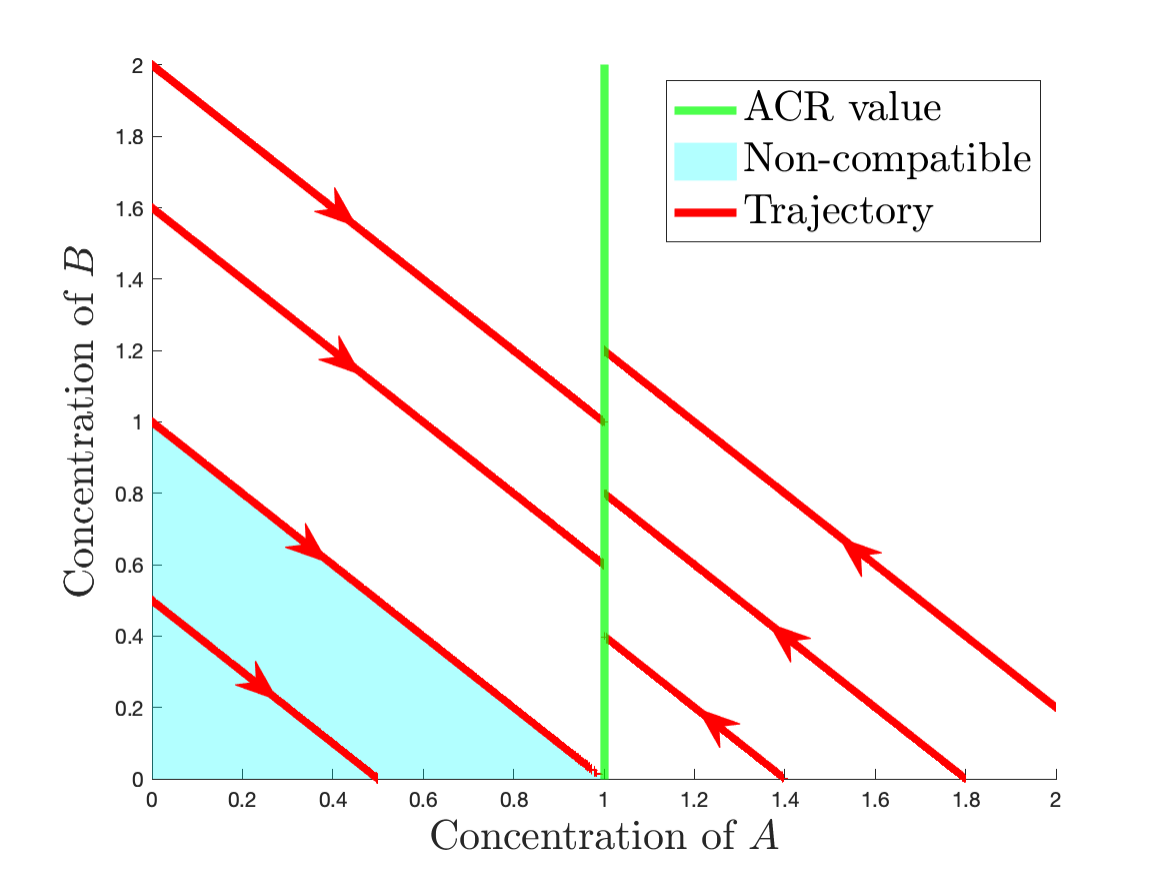}
\caption{Trajectories in phase plane}
\end{subfigure}
\caption{{\bf (Archetypal Wide Basin Dynamic ACR Network)} A dynamic ACR reaction network ($A+B \to 2B, B \to A$) with $A$ as a wide basin dynamic ACR variable. The concentration of $A$ is {\bf bounded} within the subset of $\R^2_{\ge 0}$ that is not compatible the ACR hyperplane $\{a=1\}$ (non-compatible region shown here in cyan).}
\label{fig:qperuhg305y84oigg}
\end{figure}
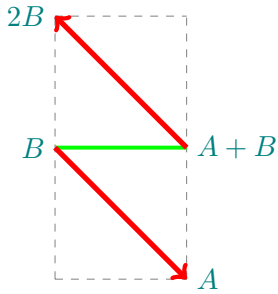
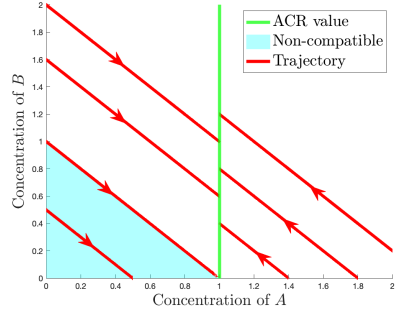

\begin{enumerate}[label={\bf Ex \arabic*}.,wide, labelwidth=!, labelindent=0pt]
\setcounter{enumi}{\value{break}}

\item\label{ex5} {\em (A minimal full basin dynamic ACR system.)}
Consider the mass action system associated to the reaction network 
\[
A+ B \stackrel[k_2]{k_1}{\rightleftarrows} B, 
\]
whose mass action ODEs are 
\[
\dot a = -k_1 ab + k_2 b, \quad \dot b = 0.
\]
Within each one-dimensional compatibility class $\{b = b(0)\}$, there is a unique globally attracting positive steady state $(a,b) = (k_2/k_1,b(0))$. So for any choice of rate constants, the resulting system is full basin dynamic ACR. The static and dynamic ACR variable is $a$ with ACR value $k_2/k_1$.

\item\label{ex6} {\em (Archetypal wide basin dynamic ACR system.)} A minimal, non-trivial, archetypal model for ACR is the network  (see also Figure \ref{fig:34p98hwrtuhw}(c) and \ref{fig:qperuhg305y84oigg}(a)): 
\begin{align*}
A + B \xrightarrow{k_1} 2B, \quad 
B \xrightarrow{k_2} A, 
\end{align*}
whose mass action ODEs are 
$
\dot a = -b(k_1 a - k_2), \quad \dot b = b(k_1 a - k_2).
$
The positive steady states form a hyperplane (ray) defined by $a=k_2/k_1$. Moreover, the positive steady states are stable and compatible with any $\{(a,b) \in \R^2_{> 0} ~|~ a+b > k_2/k_1\}$, i.e. for any $(a(0), b(0)) \in \R^2_{> 0} \setminus \{(a,b):a+b \le k_2/k_1\}$, the trajectory converges to some steady state whose $a$ coordinate is $k_2/k_1$. 
Since $\abs{\{a : a+b \le k_2/k_1\}} \le k_2/k_1$, $a$ is a wide basin dynamic ACR variable with ACR value $k_2/k_1$. See Figure \ref{fig:qperuhg305y84oigg}(b) for some sample trajectories. 

\setcounter{break}{\value{enumi}}
\end{enumerate}

\begin{figure}
\centering
\begin{subfigure}[b]{0.45\textwidth}
\begin{tikzpicture}[scale=1.75]
\draw[help lines, dashed, line width=0.25] (0,0) grid (1,2);

\node [right] at (1,2) {{\color{teal} $A+2B$}};
\node [left] at (0,1) {{\color{teal} $B$}};
\node [right] at (1,1) {{\color{teal} $A+B$}};
\node [left] at (0,0) {{\color{teal} $0$}};

\draw [-, line width=1.5, green] (0,1) -- (1,1);
\draw [->, line width=2, red] (0,1) -- (1,2);
\draw [->, line width=2, red] (1,1) -- (0,0);

\end{tikzpicture}
\caption{Reaction network embedded in Euclidean plane}
\end{subfigure}
\begin{subfigure}[b]{0.45\textwidth}
\includegraphics[scale=0.3]{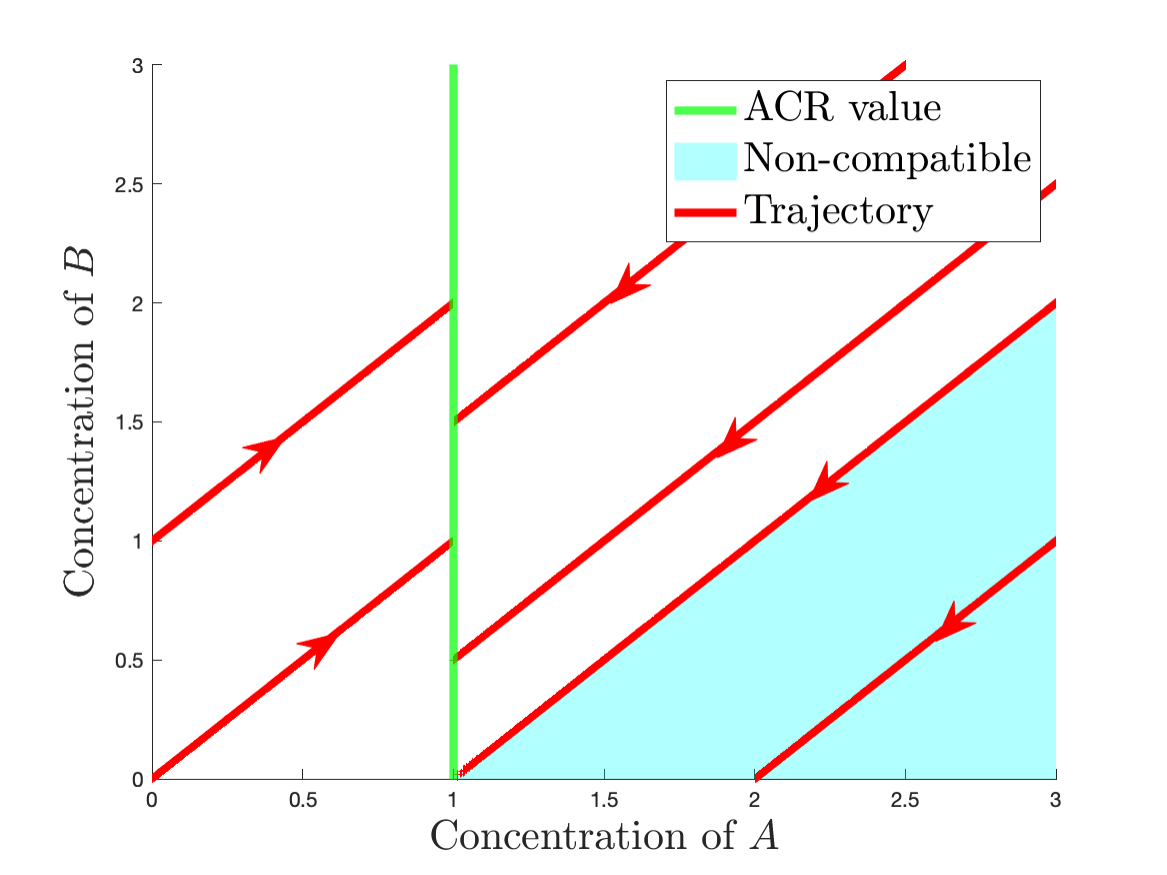}
\caption{Trajectories in phase plane}
\end{subfigure}
\caption{{\bf (Archetypal Narrow Basin Dynamic ACR Network)} A dynamic ACR reaction network ($A+B \to 0, B \to A+2B$) with $A$ as a narrow basin dynamic ACR variable. The concentration of $A$ is {\bf unbounded} within the subset of $\R^2_{\ge 0}$ that is not compatible the ACR hyperplane $\{a=1\}$ (non-compatible region shown here in cyan).}
\label{fig:pq5yhf;ghoq;35hyijg}
\end{figure}
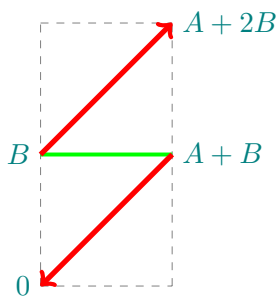
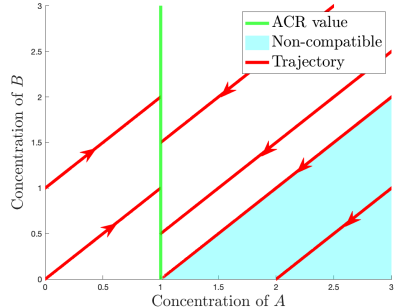

\begin{enumerate}[label={\bf Ex \arabic*}.,wide, labelwidth=!, labelindent=0pt]
\setcounter{enumi}{\value{break}}

\item\label{ex7} {\em (Static but not dynamic ACR: $f_i \equiv 0$.)}
Consider the reaction network (see also Figure \ref{fig:34p98hwrtuhw}(b))
\begin{align*}
A + B \xrightarrow{k_1} A, \quad B  \xrightarrow{k_2} 2B, 
\end{align*}
whose mass action ODE system is 
$
\dot a = 0, ~~ 
\dot b = -b(k_1 a - k_2). 
$
Here $a$ is a static ACR variable with value $k_2/k_1$. But $a$ is not a dynamic ACR variable since $\dot a \equiv 0$. 

\item\label{ex8} {\em (Static but not dynamic ACR system: Only repelling steady states.)}
Consider the reaction network (see also Figure \ref{fig:34p98hwrtuhw}(a))
\begin{align*}
A + B \xrightarrow{k_2} 2B, \quad 2A + B  \xrightarrow{k_1} 3A
\end{align*}
whose mass action ODEs are  
$
\dot a = ab(k_1 a - k_2), ~~
\dot b = -ab(k_1 a - k_2). 
$
Similar to the previous examples, the positive steady states form a hyperplane (ray) defined by $a=k_2/k_1$. The positive steady states are compatible with any $\{(a,b) \in \R^2_{\ge 0} ~|~ a+b > k_2/k_1\}$, but the steady states are unstable. The system is static ACR, with $a$ as the unique static ACR variable with ACR value $k_2/k_1$, but the system is not dynamic ACR. 

\item\label{ex9} {\em (A minimal narrow basin dynamic ACR system.)}
Consider the reaction network (see also Figure \ref{fig:pq5yhf;ghoq;35hyijg}(a) and Figure \ref{fig:24958heghuwo}(b))
\begin{align*}
A+B \xrightarrow{k_1} 0, ~~
B \xrightarrow{k_2} A + 2B,
\end{align*}
whose mass action ODEs are 
$
\dot a = -b(k_1 a - k_2), \quad 
\dot b = -b(k_1 a - k_2). 
$
The positive steady states form a ray defined by $a=k_2/k_1$. Moreover, the positive steady states are stable and compatible with any $\{(a,b) \in \R^2_{> 0} ~|~ a-b < k_2/k_1\}$. 
This shows that the system is dynamic ACR in variable $a$ with value $k_2/k_1$. 
Since $\abs{\{a: a-b \ge k_2/k_1\}}$ has no upper bound, $a$ is narrow basin dynamic ACR. See Figure \ref{fig:pq5yhf;ghoq;35hyijg}(b) for some sample trajectories. 

\setcounter{break}{\value{enumi}}
\end{enumerate}

\subsection{Higher dimensional systems.}

\begin{enumerate}[label={\bf Ex \arabic*}.,wide, labelwidth=!, labelindent=0pt]
\setcounter{enumi}{\value{break}}

\item\label{ex10} {\em (Static ACR in all variables but dynamic ACR in none.)} 
Consider the classic Lotka-Volterra system
\begin{align*}
A + B \to 2B, \quad 
B  \to 0, \quad 
A \to 2A
\end{align*}
We can apply the Shinar-Feinberg ACR criterion \cite{shinar2010structural} to this system. We check that it has deficiency $6-3-2=1$, and two non-terminal complexes differ in exactly one species. In fact, the last holds for both species, $A = (A+B) - B$ and $B = (A+B) - A$. So, by the Shinar-Feinberg ACR criterion, for all positive rate constants, the system is static ACR and concentrations of both $A$ and $B$ are static ACR variables. By Theorem \ref{thm:q035yhtihgje}, the system has a unique positive steady state for every choice of rate constants. 
The system is not dynamic ACR because for any choice of positive rate constants, the unique positive steady state is not attracting.

\item\label{ex11} {\em (Unique positive steady state which is stable for some but not all parameters: Static ACR for all rate constants; dynamic ACR system for some but not all rate constants.)} 
Consider the following reaction network, which is a simplified version of the  Sel'kov oscillator \cite{sel1968self}. 
\begin{align} \label{eq:205y8hogpo35h}
0 \xrightarrow{\rho} X \xrightarrow{\sigma} Y \xrightarrow{1} 0, \quad  
X+2Y  \xrightarrow{1} 3Y. 
\end{align}
The mass action ODE system is
\begin{align} \label{eq:4624jhwpheo}
\dot x = \rho - \sigma x - xy^2, \quad 
\dot y = -y + \sigma x + xy^2. 
\end{align}

\begin{figure}[h!]
\centering
\includegraphics[scale=0.41]{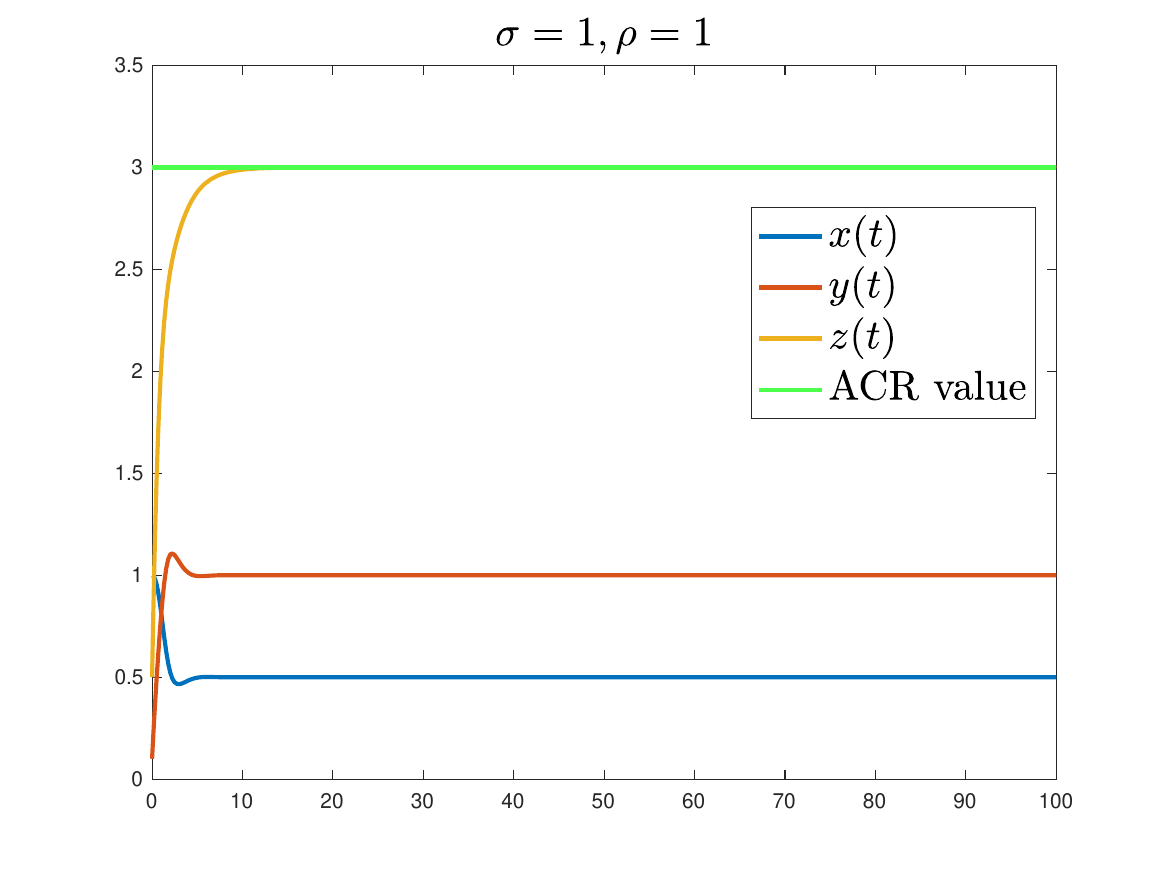}
\includegraphics[scale=0.41]{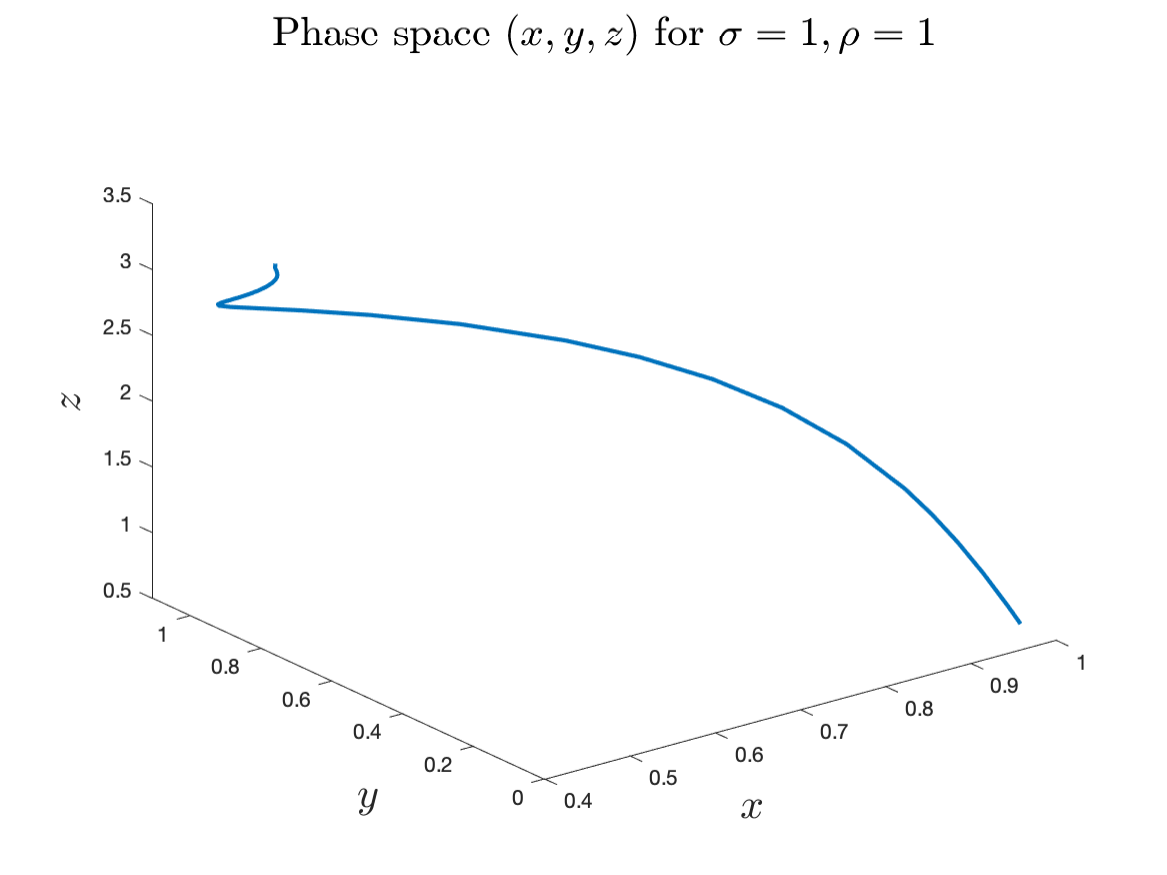}
\includegraphics[scale=0.41]{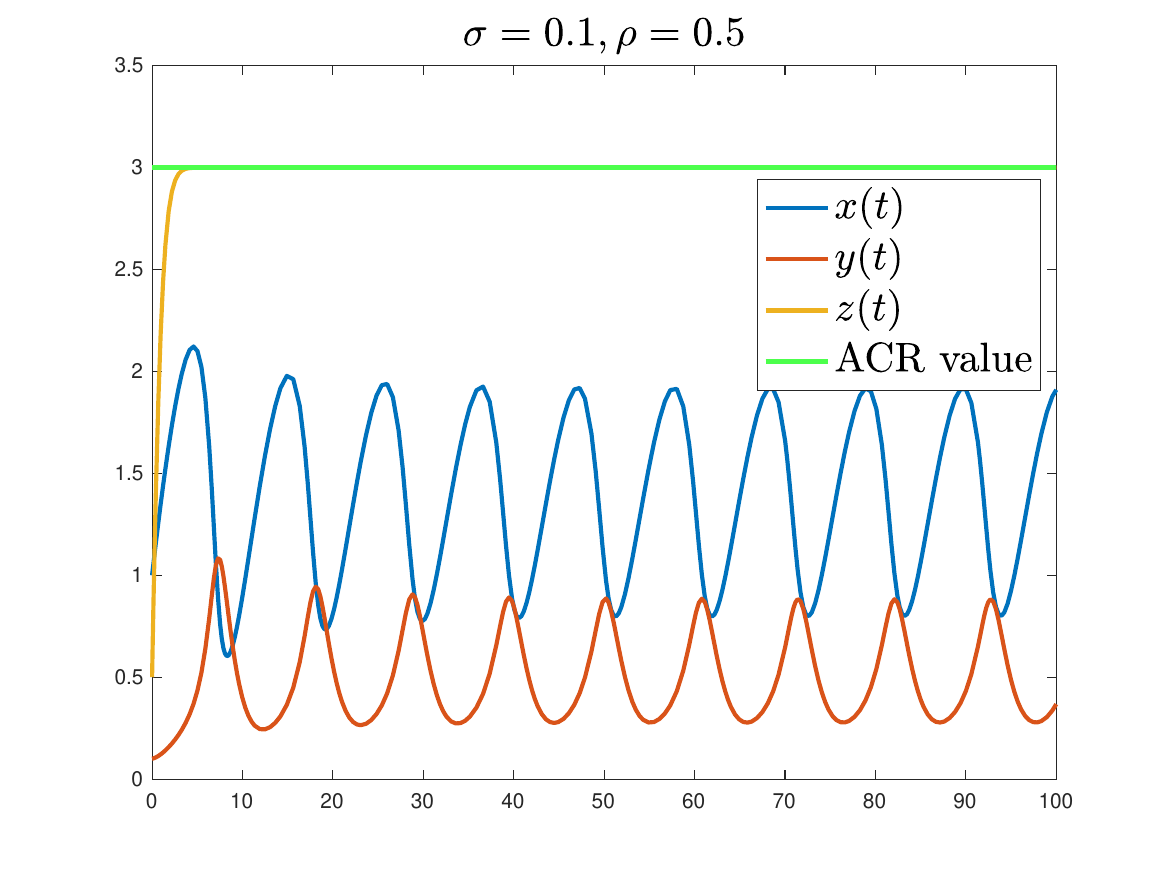}
\includegraphics[scale=0.41]{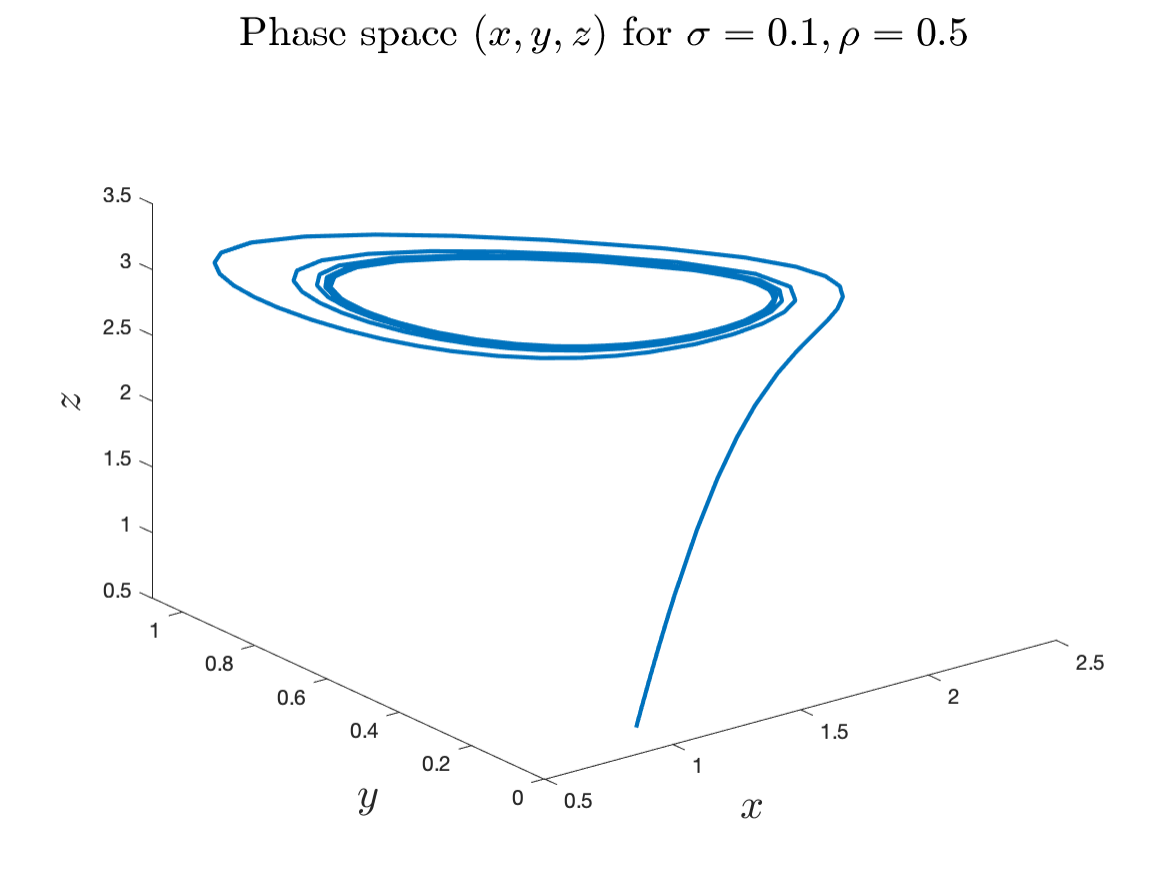}
\caption{Solutions of the system in \eqref{eq:qeoirh4toihohoewith}. {\bf (Top:)} For $\sigma =1$ and $\rho=1$, the solution converges to a positive steady state for any positive initial value. Therefore, the system is dynamic ACR, and all variables are dynamic ACR variables. {\bf (Bottom:)} For $\sigma =0.1$ and $\rho=0.5$, $z(t) \xrin k_2/k_1$, while $x(t)$ and $y(t)$ converge to non-constant periodic functions of time $t$. Therefore, the system is dynamic ACR and only $z$ is a dynamic ACR variable.}
\label{fig:67p9tkjpoyjpot}
\end{figure}

It is easy to check that for all positive rate constants, the system has a unique positive steady state whose value is $\ds \left(x^*, y^* \right) = \left(\rho/(\sigma + \rho^2), \rho \right)$. Therefore, for all positive rate constants, the system is static ACR and both $x$ and $y$ are static ACR variables. 
Furthermore, it can be checked that for $\sigma =0.1$ and $\rho=0.5$, the steady state is not attracting, so that for this choice of rate constants, the system is not dynamic ACR.  
For $\sigma =1$ and $\rho=1$, the unique positive steady state is a global attractor, and so for this choice of rate constants both $x$ and $y$ are full basin dynamic ACR variables. 

\item\label{ex12} {\em (Static ACR in all variables; dynamic ACR in one variable or all variables depending on rate constants.)} 
Consider the following reaction network, which has embedded within it the  Sel'kov oscillator.
\begin{align} \label{eq:qo35hdflfhor}
0 \xrightarrow{\rho} X \xrightarrow{\sigma} Y \xrightarrow{1} 0, \quad  
X+2Y  \xrightarrow{1} 3Y, \quad 
Z + X \stackrel[k_2]{k_1}{\rlas} X
\end{align}
The mass action ODE system is
\begin{align} \label{eq:qeoirh4toihohoewith}
\dot x = \rho - \sigma x - xy^2, \quad 
\dot y = -y + \sigma x + xy^2, \quad 
\dot z = x(k_2-k_1z)
\end{align}
It is easy to check that for all positive rate constants, the system has a unique positive steady state whose value is $\ds \left(x^*, y^*, z^* \right) = \left(\rho/(\sigma + \rho^2), \rho, k_2/k_1\right)$. Therefore, for all positive rate constants, the system is static ACR and all three variables $x,y$, and $z$ are static ACR variables.  
Moreover, for any choice of positive rate constants and for any initial value, $z \xrin k_2/k_1$. 
It follows then that for any choice of positive rate constants, the system is dynamic ACR, and that $z$ is a full basin dynamic ACR variable whose dynamic ACR value is $k_2/k_1$. Note that this value is the same as the static ACR value of $z$. 
Furthermore, it can be checked that for $\sigma =0.1$ and $\rho=0.5$, the steady state is not attracting, so that for this choice of rate constants, the only dynamic ACR variable is $z$.  
For $\sigma =1$ and $\rho=1$, the unique positive steady state is a global attractor, and so for this choice of rate constants all three variables, $x, y$ and $z$ are full basin dynamic ACR variables. See Figure \ref{fig:67p9tkjpoyjpot} for the trajectories for the two choices of rate constants. 

\setcounter{break}{\value{enumi}}
\end{enumerate}

\subsection{Dynamic ACR but not static ACR.}

\begin{enumerate}[label={\bf Ex \arabic*}.,wide, labelwidth=!, labelindent=0pt]
\setcounter{enumi}{\value{break}}
 
\item\label{ex401} {\em (Adding reactions with unrelated species can destroy static ACR but always preserves dynamic ACR.)} 
The system $0 \stackrel[k']{k}{\rlas} A$ in \ref{ex4}, is both static and dynamic ACR for any choice of rate constants. Suppose we add a flow reaction of the type $0 \xrightarrow{g} B$, so that the new mass action system is $\dot a = k - k'a, ~\dot b = g >0$. 
Then since the concentration of $B$ goes to infinity there are no positive steady states, and so the system is not static ACR. However $a$ is still a full basin dynamic ACR variable with the same value $k/k'$. 
Similar considerations apply if instead of the inflow $0 \xrightarrow{g} B$, we add the outflow reaction $B \xrightarrow{\ell} 0$. 
\setcounter{break}{\value{enumi}}
\end{enumerate}
The result in the previous example holds in general. 
\begin{proposition} \label{prop:simpleunion}
Consider two dynamical systems $\DD_x$ given by $\dot x = f(x)$ with $x \in \R_{\ge 0}^n$ and $\DD_y$ given by $\dot y = g(y)$ with $y \in \R_{\ge 0}^m$. Suppose that $\DD_x$ has a dynamic ACR variable $x_i$ with value $a_i^*$.  Then the dynamical system $\DD_x \cup \DD_y$ given by $\{\dot x = f(x), ~\dot y = g(y)\}$ with $(x,y) \in \R_{\ge 0}^{n+m}$ has $x_i$ as a dynamic variable with the same value $a_i^*$. 
\end{proposition}
Even when the dynamics are bounded and do not converge to the boundary, it is possible to have a dynamic ACR system which is not static ACR, as the following example shows. 

\begin{enumerate}[label={\bf Ex \arabic*}.,wide, labelwidth=!, labelindent=0pt]
\setcounter{enumi}{\value{break}}

\item\label{ex13} 
{\em (Dynamic ACR but not static ACR in a mass conserving system.)} 
Consider the mass action dynamical system resulting from the following reaction network, where the labels on the arrows indicate reaction rate constants. 
\begin{align} \label{eq:peogueghigroe}
A + 2B &\xrightarrow{2} 2A + B \xrightarrow{1} 3A \nonumber \\ 
3A + B &\xrightarrow{1} 2A + 2B \xrightarrow{2} A + 3B 
\end{align}
The resulting system of ODEs is:
\begin{align*}
\dot a &= -ab(a+2b)(a-1) \\ 
\dot b &= ab(a+2b)(a-1)
\end{align*}
Note that $(a + b)(t)$ is a constant function of time $t$.
The variable $a$ is static and full basin dynamic ACR with the ACR value of $1$, since clearly for any positive initial condition $a \xrin 1$, and $a=1$ gives a unique positive steady state that is compatible with any initial condition with $a(0) + b(0) > 1$. 

Now consider the same reaction network as \eqref{eq:peogueghigroe} with one additional reaction, and the reaction rate constants as shown below. 
\begin{align} \label{eq:eoirhhgoierhg}
&A + B \xrightarrow{1} 2B \nonumber \\
A + 2B &\xrightarrow{2} 2A + B \xrightarrow{2} 3A \nonumber \\ 
3A + B &\xrightarrow{1} 2A + 2B \xrightarrow{2} A + 3B 
\end{align}
The resulting system of ODEs is:
\begin{align*}
\dot a &= -ab(a+2b-1)(a-1) \\ 
\dot b &= ab(a+2b-1)(a-1)
\end{align*}
Once again, $(a + b)(t)$ is a constant function of time $t$. The system has positive steady states given by $a=1$ or $a+2b=1$ (see Figure \ref{fig:eorigheq5h095039} (left)). Clearly, then the system is not static ACR. 
For a positive initial condition $(a(0), b(0))$ to be compatible with $a=1$, it must be the case that $a(0) + b(0) > 1$ which implies that $a(t) + 2b(t) \ge a(t) + b(t) = a(0) + b(0)  > 1$, and so again for any positive initial condition compatible with $\{a=1\}$, we have $a \xrin 1$. 
Therefore the system is dynamic ACR with $a$ as a wide basin dynamic ACR variable with ACR value $1$. 

\begin{figure}
\centering
\includegraphics[scale=0.41]{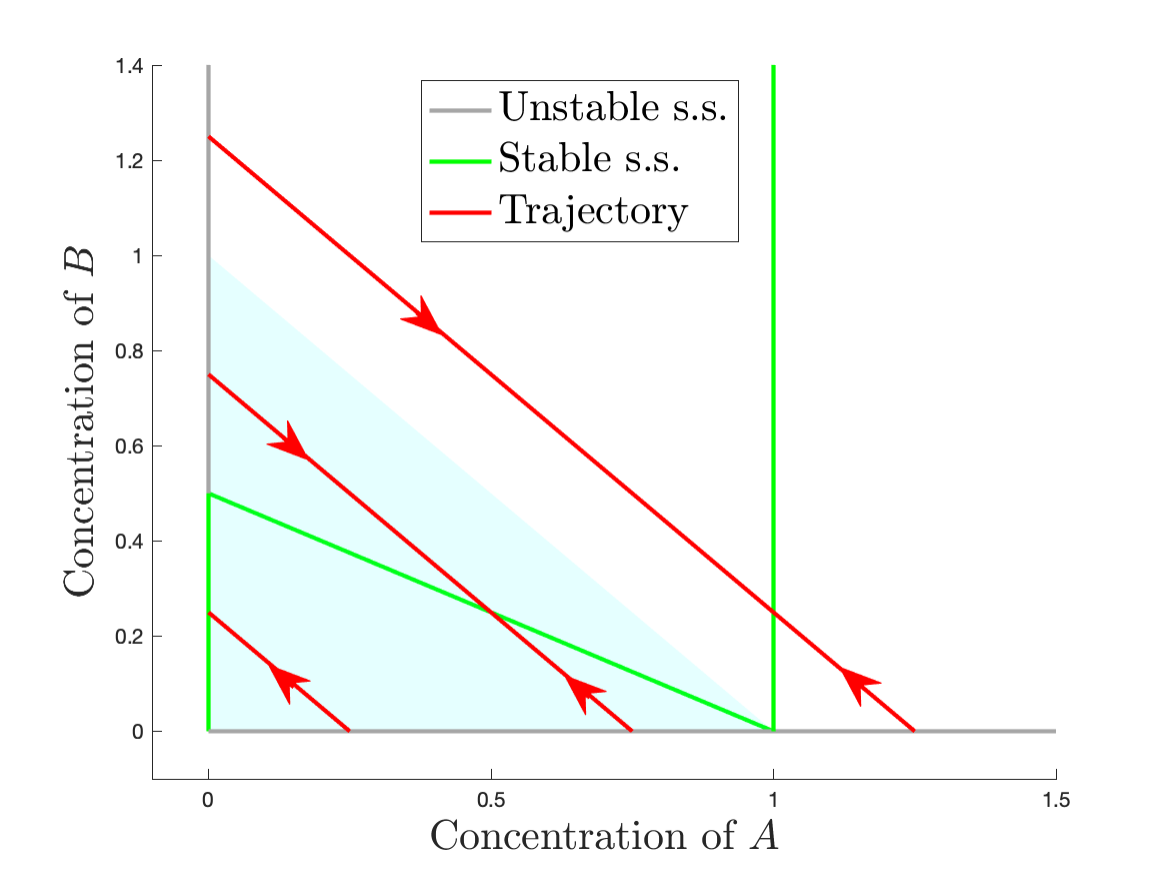}
\includegraphics[scale=0.41]{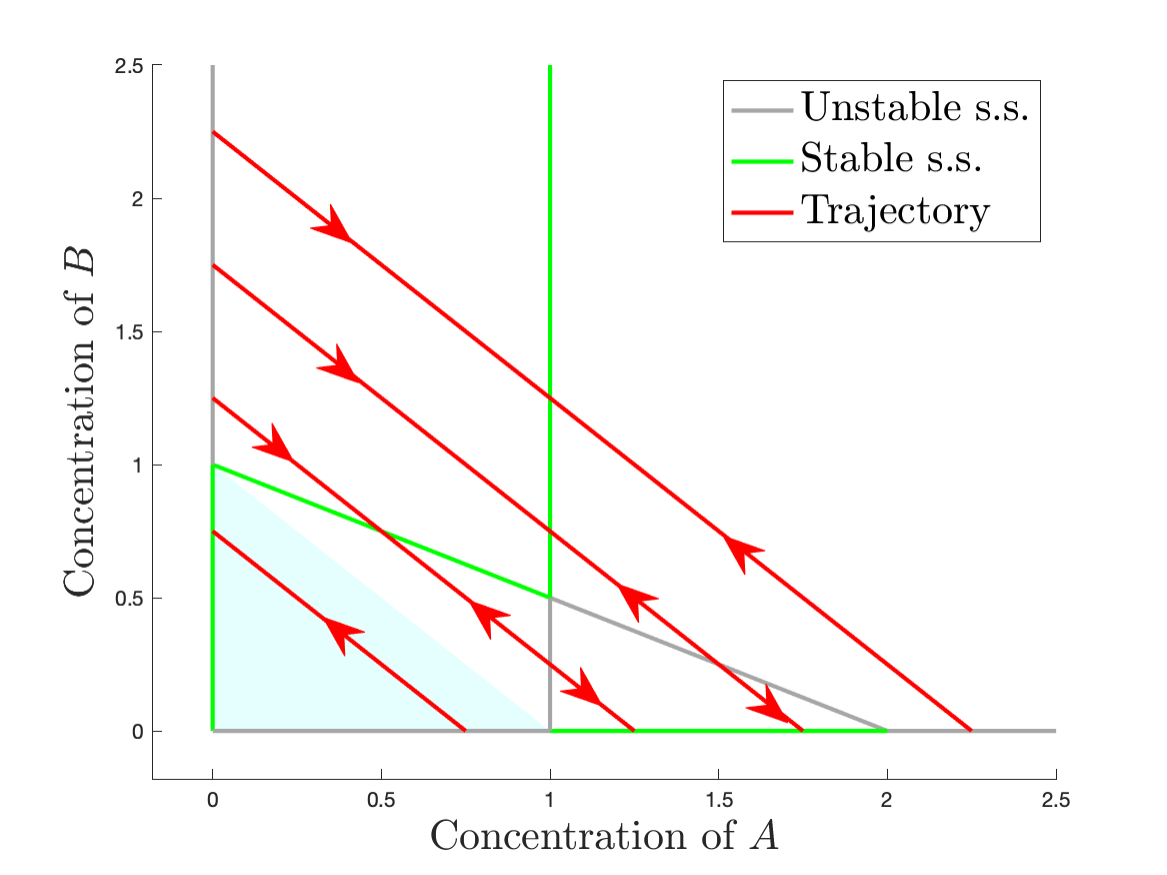}
\caption{Systems that are not static ACR because there are additional steady states outside a positive hyperplane (ray). {\bf (Left:)} Dynamic ACR system. All positive initial values compatible with $\{a=1\}$ result in convergence  to $a=1$. {\bf (Right:)} Not dynamic ACR. For sufficiently large positive initial conditions (in this case $a(0) + b(0) >2$), there is convergence to $a=1$.}
\label{fig:eorigheq5h095039}
\end{figure}

Finally, consider the same reaction network as \eqref{eq:eoirhhgoierhg} but with different rate constants as below. 
\begin{align} \label{eq:erogiheoghdsg}
&A + B \xrightarrow{2} 2B \nonumber \\
A + 2B &\xrightarrow{2} 2A + B \xrightarrow{3} 3A \nonumber \\ 
3A + B &\xrightarrow{1} 2A + 2B \xrightarrow{2} A + 3B 
\end{align}
The resulting system of ODEs is:
\begin{align*}
\dot a &= -ab(a+2b-2)(a-1) \\ 
\dot b &= ab(a+2b-2)(a-1)
\end{align*}
This system is again not static ACR. Furthermore, if $1 < a(0) + b(0) < 2$, then $(a(0),b(0))$ is compatible with $\{(a,b) \in \R^2_{>0} ~|~ a=1\}$, but $a \not \to 1$ for every such initial condition. Finally, $a \xrin 1$ if $a(0) + b(0) > 2$, and thus we do have convergence to a steady state with $a=1$ outside a compact set (see Figure \ref{fig:eorigheq5h095039} (right)). 

\setcounter{break}{\value{enumi}}
\end{enumerate}

\begin{remark}
The last example suggests a naturally arising weaker version of dynamic ACR. For instance, we might want to say that \eqref{eq:erogiheoghdsg} is dynamic ACR on $\{a+b>2\}$.  We study this and other weaker forms of dynamic ACR in \cite{joshi2022motifs}. 

If a small positive inflow parameter is added either to $\dot a$ or $\dot b$ (corresponding to influx of species $A$ or $B$ into the system), one can show that the resulting system is dynamic ACR in all three cases \eqref{eq:peogueghigroe}, \eqref{eq:eoirhhgoierhg}, and \eqref{eq:erogiheoghdsg}. 
\end{remark}

\subsection{Wide basin dynamic ACR.}
The condition for a dynamic ACR variable to be wide basin dynamic ACR is that within the set of positive points that are incompatible with the steady states, the ACR variable is bounded. 
In other words, if the initial concentration of the ACR variable is sufficiently large, the ACR variable will converge to its ACR value. 
Even so, the set of incompatible points can be both empty or unbounded. We start with an example of an unbounded case.

\begin{enumerate}[label={\bf Ex \arabic*}.,wide, labelwidth=!, labelindent=0pt]
\setcounter{enumi}{\value{break}}

\item\label{ex14} {\em (Set of initial values that do not converge to ACR value may be unbounded.)}
The ACR variable has an upper bound in the set of initial values that do not converge to the ACR value, by definition. But the other variables may not have any upper bound as the following example shows. 
Consider the reaction network:
\begin{align*}
X + Y + Z &\xrightarrow{k_1} 2X + 2Y \\
X + Y &\xrightarrow{k_2} Z
\end{align*}
The stoichiometric subspace is the span of the vector $(1,1,-1)$, so there are two independent conservation relations, for instance $x+z = c_1$ and $y+z=c_2$. The system of ODEs is 
\[
-\dot x = -\dot y = \dot z = -k_1xy(z-k^*), 
\]
where $k^* = k_2/k_1$. It's easy to see that $z$ is a dynamic ACR variable with ACR value $k^*$. Furthermore, the following result holds, which we state without proof since the proof is easy. 
\begin{claim} 
Consider the system $-\dot x = -\dot y = \dot z = -k_1xy(z-k^*)$, with $c_1 = x(0) + z(0) > 0$ and $c_2 = y(0) + z(0) > 0$. Let $c \coloneqq \min\{c_1,c_2\}$. If $c \le k^*$, then $(x(t),y(t),z(t)) \xrin (c_1-c, c_2-c, c)$, i.e. convergence is to a boundary steady state. If $c > k^*$, then $(x(t),y(t),z(t)) \xrin (c_1-k^*, c_2-k^*, k^*)$, i.e. convergence is to a positive steady state with $z$ converging to the ACR value. 
\end{claim}
\begin{proof}
The solution of the initial value problem $\DD = \{-\dot x = -\dot y = \dot z = -k_1xy(z-k^*), x(0) = x_0, y(0) = y_0, z(0) = z_0\}$ and that of $\wt \DD = \{-\dot{\wt x} = -\dot{\wt y} = \dot{\wt z} = -(\wt z-k^*), \wt x(0) = x_0, \wt y(0) = y_0, \wt z(0) = z_0\}$ are equivalent as trajectories when restricted to the nonnegative orthant $\R^3_{\ge 0}$. 
The initial value problem $\wt \DD$ can be explicitly solved: 
\begin{align*}
\wt x(t) &= \left(x_0 + z_0 - k^*\right) - \left(z_0 - k^* \right) e^{-t} \\
\wt y(t) &= \left(y_0 + z_0 - k^*\right) - \left(z_0 - k^* \right) e^{-t} \\
\wt z(t) &= k^* + \left(z_0 - k^* \right) e^{-t}.
\end{align*}
If $c > k^*$, then the set of points $\{(\wt x(t), \wt y(t), \wt z(t)) :{t \ge 0}\}$ is entirely contained in the positive orthant $\R^3_{> 0}$ and so $\{(\wt x(t), \wt y(t), \wt z(t)) :{t \ge 0}\} = \{(x(s), y(s), z(s)) :{s \ge 0}\}$ as a set and moreover, $\lim_{s \to \infty} (x(s), y(s), z(s)) = \lim_{t \to \infty} (\wt x(t), \wt y(t), \wt z(t)) = (c_1-k^*, c_2-k^*, k^*)$. 

On the other hand, if $c \le k^*$, then  $\lim_{t \to \infty} (\wt x(t), \wt y(t), \wt z(t)) = (c_1-k^*, c_2-k^*, k^*) \notin \R^3_{> 0}$. 
Therefore, $\lim_{s \to \infty} (x(s), y(s), z(s))$ is either the intersection point of $\{(\wt x(t), \wt y(t), \wt z(t)) :{t \ge 0}\}$ with the boundary of $\R^3_{\ge 0}$ (when $c < k^*$) or equal to $\lim_{t \to \infty} (\wt x(t), \wt y(t), \wt z(t)) = (c_1-k^*, c_2-k^*, k^*)$ (when $c = k^*$). In either case, $(x(t),y(t),z(t)) \xrin (c_1-c, c_2-c, c)$. 
\end{proof}

Notice in particular the region of initial conditions that  do not result in convergence to the ACR value is $\min\{c_1,c_2\} \le k^*$ which includes points with small $x$ and $z$ values but arbitrarily large $y$ values. Thus the set $S = \{s(0) \in \R^3_{\ge 0} ~|~ s(0) \mbox{ not compatible with } \{z = k^*\} \cap \R^3_{\ge 0} \}$ is non-compact. However $\{z(0) ~|~ (x(0),y(0),z(0)) \in S\}$ is bounded above by $k^*$, which shows that $z$ is a wide basin dynamic ACR variable. 

\item\label{ex15} {\em (Set of initial values that do not converge to ACR may be empty: full basin dynamic ACR.)}
The flip side of the previous category is the class of full basin dynamic ACR systems, for which every positive initial condition converges to an ACR value. We consider a simplified model of the ground state of a carbon nanotube rope \cite{cobden1998spin}. The ground state may have fractional or integral spin, alternating between the two as a spin $+1/2$ electron is absorbed or emitted by the carbon nanotube rope. 

Suppose that $c_e$  and $c_o$ denote the concentration of carbon nanotube ropes that have even and odd number of electrons, respectively, in the ground state. Suppose $x$ is the concentration of free electrons in the ambient space. We assume that electrons are absorbed or emitted with a rate constant that depends only on the odd or even state. Then we can represent the system as a reaction network. 
\begin{align*}
X+C_e \stackrel[k_2]{k_1}{\rightleftarrows} C_o, \quad X+C_o \stackrel[k_4]{k_3}{\rightleftarrows} C_e
\end{align*}
The system of mass action ODEs associated with the network is:
\begin{align*}
\dot x &= -k_1xc_e + k_2 c_o - k_3xc_o + k_4c_e,\\ 
\dot c_e &= -k_1xc_e + k_2 c_o + k_3xc_o - k_4c_e,\\ 
\dot c_o &= k_1xc_e - k_2 c_o - k_3xc_o + k_4c_e.
\end{align*}
The quantity $c_e+c_o$ is conserved over time, so the dynamics are restricted to a two-dimensional affine set or compatibility class $S_c \coloneqq \{(x,c_e,c_o) \in \R^3_{\ge 0} : c_e + c_o = c > 0 \}$. The deficiency is $4-2-2=0$ and the network is weakly reversible (i.e. every reaction is part of a cycle). Many things are known about the dynamical properties of reversible, deficiency $0$ systems \cite{horn1972necessary,feinberg2019foundations}, see also Section \ref{sec:5yiu6fikfyuerlkgfye}. For instance, reversible, zero deficiency systems have a unique positive steady state within each compatibility class $\{c>0\}$, and each of these steady states attracts all compatible, positive initial values. To explicitly solve for the positive steady state, note that 
each of the two binomials $k_1xc_e - k_2 c_o$ and $k_3xc_o - k_4c_e$ must vanish at the steady state. So, if we denote a steady state by $(x^*, c_e^*, c_o^*)$ then we have 
\begin{align*}
x^* = \frac{k_2 c_o^*}{k_1c_e^*} = \frac{k_4 c_e^*}{k_3c_o^*} \implies \frac{c_e^*}{c_o^*} = \sqrt{\frac{k_2k_3}{k_1k_4}} \implies x^* = \sqrt{\frac{k_2k_4}{k_1k_3}}. 
\end{align*}
Note that the value of $x^*$ is independent of $c$. From general results about reversible, zero deficiency systems, $x$ is a dynamic ACR variable with the dynamic ACR value $\sqrt{k_2k_4/(k_1k_3)}$. 
Thus the basin of attraction of $\{x^* = \sqrt{k_2k_4/(k_1k_3)}\}$ is the entire positive orthant. 
Note further that while neither $c_e$ nor $c_o$ is an ACR variable, their ratio $c_e/c_o$ behaves as an ACR variable. This type of ``ratio ACR'' will be discussed in detail in future work. 
\setcounter{break}{\value{enumi}}
\end{enumerate}

\section{Static and Dynamic ACR reaction networks} \label{sec:3pyihoeqhljd}

\begin{table}
\setlength\extrarowheight{0pt}
\centering
\begin{adjustbox}{max width=\textwidth}
 \begin{tabular}{|c|c||c|c||c||c||c||c|}
 \hline
 \specialcell{Capacity for\\static ACR?}& \specialcell{Is network\\static ACR?}&   \specialcell{Capacity for\\dynamic ACR?} & \specialcell{Is network\\dynamic ACR?} & \specialcell{Network} &\specialcell{Static\\ACR\\species} & \specialcell{Dynamic\\ACR\\species} & \specialcell{Location}   \\
 \hline
 \hline
Yes & Yes &Yes & Yes & $0 \rlas A$ & $A$ & $A$   &   \ref{ex4}   \\
 \hline
Yes & Yes &Yes & No & \specialcell{$0 \to A \to B \to 0$\\$A + 2B \to 3B$} & $A,B$ & $~$   & \ref{ex11}  \\
 \hline
Yes & Yes &No & No & \specialcell{$0 \la A$\\$2A \to 3A$} & $A$ & $~$   &   \ref{ex3}  \\
 \hline
 Yes & No &Yes & Yes & \specialcell{$0 \rlas A$\\$A+B \to A$\\$B \to 2B$} & $~$ & $A$   &   \ref{ex17}  \\
 \hline
 Yes & No &Yes & No & \specialcell{$0 \rlas A$\\$2A \rlas 3A$}  & $~$ & $~$  &   \ref{ex102}  \\
 \hline 
 Yes & No & No & No & \specialcell{$0 \rlas A$\\$2A \to 3A$}  & $~$ & $~$  &   \ref{ex101}  \\
 \hline 
No & No & Yes & Yes & \specialcell{$0 \to A+B$\\$A \to B$} & $~$ & $A$   &   \ref{ex201}  \\
 \hline 
 No & No & Yes & No & \specialcell{$0 \to A+B$\\$A \to B$\\$A \to 2A$} & $~$ & $~$   &   \ref{ex202}  \\
 \hline 
  No & No & No & No & \specialcell{$A \rlas B$} & $~$ & $~$   &   \ref{ex203}  \\
 \hline 
 \end{tabular}
 \end{adjustbox}
 \caption{{\em (All possibilities are realized)} For a given network $\GG$, the set of questions: (1) does $\GG$ have capacity for static ACR?, (2) is $\GG$ static ACR? has 3 distinct possible answers: (yes, yes), (yes, no) and (no, no), while (no, yes) is ruled out since a static ACR network must have the capacity for static ACR by definition. Similarly there are 3 distinct answers when the word ``static'' is replaced by ``dynamic'' in the previous sentence. Combining these two observations, the four questions in the first four rows have 9 distinct answers -- {\em all of which are realized in reaction networks}, as this table shows. Furthermore, the examples are relatively simple since they involve only one or two species and very few reactions.} \label{tab:58yhoithosigjpi3}
 \end{table}

We are often interested in mass action systems resulting from reaction networks. It is possible that the mass action system resulting from a reaction network is a dynamic ACR system for one choice of rate constants but not for another choice. Moreover, it may be that a specific variable is dynamic ACR for only a proper subset of choices of rate constants for which the overall system is dynamic ACR. Similar remarks may hold even for static ACR. For instance, a bifurcation from a unique positive steady state to two or more steady states is likely to turn a static ACR system into one that is not static ACR. 
It behooves us to single out for special attention a reaction network that results in a static or dynamic ACR system for any choice of rate constants.

\begin{definition} \label{def:049jhiwr}
Suppose that $(\GG,K)$ is a mass action system resulting from the reaction network $\GG$, where $K$ denotes the specific choice of mass action rate constants. 
\beit
\item We say that $\GG$ has capacity for static (dynamic) ACR if there is a $K$ such that the mass action system $(\GG,K)$ is a static (dynamic) ACR system.
\item We say that $\GG$ is a {\em static (dynamic) ACR network} if $(\GG,K)$ is a static (dynamic) ACR system for all choices of $K$. 
\item We say a species $X$ in a network $\GG$ is a {\em static (respectively: dynamic, wide basin dynamic, narrow basin dynamic, full basin dynamic) ACR species} if the concentration of $X$ is a static (respectively: dynamic, wide basin dynamic, narrow basin dynamic, full basin dynamic) ACR variable in $(\GG,K)$ for all choices of $K$. 
\enit
\end{definition}
\begin{remark}
The above definition deviates somewhat from the Shinar-Feinberg definition of an ACR species. Since they restrict attention to a fixed mass action system (i.e. fixed rate constants), their notion of an ACR species is analogous to our notion of an ACR variable in a dynamical system. 
\end{remark}

We now set about to present examples of networks which answer to these questions: (i) does the network have capacity for static ACR? (ii) is the network static ACR? (iii) does the network have capacity for dynamic ACR? (iv) is the network dynamic ACR? Of the $2^4=16$ distinct possibilities of yes or no answers to these questions, we can {\it a priori} rule out `no' to (i) and `yes' to (ii) -- a static ACR network must have capacity for static ACR. We can similarly rule out `no' to (iii) and `yes' to (iv). We are left then with $3 \times 3 = 9$ distinct possibilities of yes/no answers to the four questions. Each of these possibilities is realized in a fairly simple network involving no more than 2 species and no more than 4 reactions. The results are presented in Table \ref{tab:58yhoithosigjpi3}.

\begin{enumerate}[label={\bf Ex \arabic*}.,wide, labelwidth=!, labelindent=0pt]
\setcounter{enumi}{\value{break}}

\item\label{ex101} {\em (Has capacity for static ACR but not for dynamic ACR)}
Consider the reaction network
\begin{align} \label{eq:itjoihjoirwtej}
0 \stackrel[k_2]{k_1}{\rlas} A, \quad 2A \xrightarrow{k_3} 3A. 
\end{align}
The mass action ODE system is
\begin{align} \label{eq:w4pyjdligjoiej}
\dot a=k_1-k_2a+k_3a^2. 
\end{align}
The mass action system has a unique positive steady state if and only if the rate constants satisfy $k_2^2=4k_1k_3$. 
Furthermore, this unique steady state is unstable. This implies that when the constraint is satisfied, the system is static ACR but not dynamic ACR.
For any other choice, there are either no positive steady states or two distinct positive steady states. In either case, the system is neither static ACR nor dynamic ACR. 
We conclude that the reaction network has capacity for static ACR, is not a static ACR reaction network, and has no capacity for dynamic ACR.

\item\label{ex102} {\em (Has capacity for static and dynamic ACR but neither static nor dynamic ACR network due to multistationarity)}
Consider the reaction network
\begin{align} \label{eq:4o5yohegodo}
0 \stackrel[k_2]{k_1}{\rlas} A, \quad 2A \stackrel[k_4]{k_3}{\rlas} 3A. 
\end{align}
The mass action ODE system is
\begin{align} \label{eq:35pp7ulkdkuhe}
\dot a=k_1-k_2a+k_3a^2-k_4a^3. 
\end{align}
It is fairly easy to show that for $k_i=1$ for $i \in \{1,2,3,4\}$, there is a unique positive steady state which is globally attracting. So the system has capacity for static and dynamic ACR, and the variable $a$ is full basin dynamic ACR for this choice of rate constants. 

For $k_1=4, k_2=8, k_3=3.5, k_4=0.4$, there are three positive steady states and so the resulting system is neither static nor dynamic ACR. Therefore, the network is neither static nor dynamic ACR.
\setcounter{break}{\value{enumi}}
\end{enumerate}
A network that satisfies the Shinar-Feinberg criterion may lack steady states for all choices of rate constants. Such a network obviously has no capacity for static ACR. Somewhat more surprisingly, there exist networks that satisfy the Shinar-Feinberg criterion and have steady states for some but not all choices of rate constants. Such a network has the capacity for static ACR but is not static ACR. We present an example below. 
\begin{enumerate}[label={\bf Ex \arabic*}.,wide, labelwidth=!, labelindent=0pt]
\setcounter{enumi}{\value{break}}
\item\label{ex16} {\em (Has capacity for static and dynamic ACR but neither static nor dynamic ACR network due to possible absence of steady states)}
Consider the reaction network
\begin{align} \label{eq:ioyuktjhp}
0 \xleftarrow{k_1} A \stackrel[k_3]{k_2}{\rlas} 2A 
\end{align}
The mass action ODE system is
\begin{align} \label{eq:4o68yigklpw}
\dot a = a\left(k_2 - k_1 -k_3a \right)
\end{align}
If $k_2 \le k_1$ then there are no positive steady states and the system is neither static nor dynamic ACR. If $k_2>k_1$, then $a^*=(k_2-k_1)/k_3$ is the unique positive steady state which means that $a$ is a static ACR variable with ACR value $(k_2-k_1)/k_3$. In fact, in this case $a$ is a full basin dynamic ACR variable.

\item\label{ex201} {\em (Has no capacity for static ACR but has capacity for dynamic ACR)}
Consider the reaction network
\begin{align} \label{eq:ehrguoehgld}
0 \xrightarrow{k_1} A+B, \quad  A \xrightarrow{k_2} B, \quad  A \xrightarrow{k_3} 2A.   
\end{align}
The mass action ODE system is
\begin{align} \label{eq:24o58yjiegr}
\dot a = k_1 - \left(k_2 - k_3\right) a, \quad 
\dot b = k_1 + k_2 a
\end{align}
Clearly $b \xrin \infty$, so the system has no capacity for static ACR. If $k_2 \le k_3$, then $a \xrin \infty$, so in this case the system is not dynamic ACR. But if $k_2 > k_3$, then $a \xrin k_1/(k_2-k_3)$ for any positive initial condition, which implies that the system is full basin dynamic ACR. 

\item\label{ex202} {\em (Has no capacity for static ACR but is dynamic ACR)}
Consider the reaction network
\begin{align} \label{eq[42iyjjfskgjejg}
0 \xrightarrow{k_1} A+B, \quad  A \xrightarrow{k_2} B.   
\end{align}
The mass action ODE system is
\begin{align} \label{eq:oeprighoqehgr}
\dot a = k_1 - k_2 a, \quad 
\dot b = k_1 + k_2 a
\end{align}
Clearly $b \xrin \infty$, so the system has no capacity for static ACR. On the other hand, $a \xrin k_1/k_2$ for any positive initial condition, which implies that the system is full basin dynamic ACR for any choice of rate constants. Therefore, $A$ is a full basin dynamic ACR species. 
\setcounter{break}{\value{enumi}}
\end{enumerate}

\begin{figure}[h!] 
\centering
\begin{subfigure}[b]{0.2\textwidth}
\begin{tikzpicture}[scale=1]
\draw[help lines, dashed, line width=0.25] (0,0) grid (2,2);
\draw [<->, line width=1.5, blue] (0,2.5) -- (0,0) -- (2.5,0);

\node [below] at (2,0) {{\cbl $A$}};
\node [left] at (0,2) {{\cbl $B$}};

\draw [->, line width=2, red] (0,1) -- (0,2);
\draw [->, line width=2, red] (1,1) -- (1,0);

\draw [-, line width=2, green] (0,1) -- (1,1);

\end{tikzpicture}
\caption{\bf Static ACR {\large \bf ({\color{green} \cmark})} \\Dynamic ACR {\large \bf ({\color{red} \xmark})} \\Wide basin ACR {\large \bf ({\color{red} \xmark})}\\Full basin ACR {\large \bf ({\color{red} \xmark})}}
\end{subfigure}
\begin{subfigure}[b]{0.2\textwidth}
\begin{tikzpicture}[scale=1]
\draw[help lines, dashed, line width=0.25] (0,0) grid (2,2);
\draw [<->, line width=1.5, blue] (0,2.5) -- (0,0) -- (2.5,0);

\node [below] at (2,0) {{\cbl $A$}};
\node [left] at (0,2) {{\cbl $B$}};

\draw [->, line width=2, red] (1,1) -- (0,0);
\draw [->, line width=2, red] (0,1) -- (1,2);

\draw [-, line width=2, green] (0,1) -- (1,1);

\end{tikzpicture}
\caption{\bf Static ACR {\large \bf ({\color{green} \cmark})} \\  Dynamic ACR {\large \bf ({\color{green} \cmark})} \\Wide basin ACR {\large \bf ({\color{red} \xmark})}\\Full basin ACR {\large \bf ({\color{red} \xmark})}}
\end{subfigure}
  \begin{subfigure}[b]{0.2\textwidth}
\begin{tikzpicture}[scale=1]
\draw[help lines, dashed, line width=0.25] (0,0) grid (2,2);
\draw [<->, line width=1.5, blue] (0,2.5) -- (0,0) -- (2.5,0);

\node [below] at (2,0) {{\cbl $A$}};
\node [left] at (0,2) {{\cbl $B$}};

\draw [->, line width=2, red] (1,1) -- (0,2);
\draw [->, line width=2, red] (0,1) -- (1,0);

\draw [-, line width=2, green] (0,1) -- (1,1);

\end{tikzpicture}
\caption{\bf Static ACR {\large \bf ({\color{green} \cmark})} \\  Dynamic ACR {\large \bf ({\color{green} \cmark})} \\Wide basin ACR {\large \bf ({\color{green} \cmark})}\\Full basin ACR {\large \bf ({\color{red} \xmark})}}
\end{subfigure}
\begin{subfigure}[b]{0.2\textwidth}
\begin{tikzpicture}[scale=1]
\draw[help lines, dashed, line width=0.25] (0,0) grid (2,2);
\draw [<->, line width=1.5, blue] (0,2.5) -- (0,0) -- (2.5,0);

\node [below] at (2,0) {{\cbl $A$}};
\node [left] at (0,2) {{\cbl $B$}};

\draw [->, line width=2, red] (1,1) -- (0,1);
\draw [->, line width=2, red] (0,1) -- (1,1);

\end{tikzpicture}
\caption{\bf Static ACR {\large \bf ({\color{green} \cmark})} \\  Dynamic ACR {\large \bf ({\color{green} \cmark})} \\Wide basin ACR {\large \bf ({\color{green} \cmark})}\\Full basin ACR {\large \bf ({\color{green} \cmark})}}
\end{subfigure}
\caption{Minimal (in terms of total stoichiometry), non-trivial (has more than one positive steady state) static ACR networks.  $A$ is a static ACR species for each network and also a dynamic ACR species for networks (b)-(d). 
The reactant polytope for each network is the  line segment joining the complexes $B$ and $A+B$.
The total stoichiometry for each network is 6.}
\label{fig:24958heghuwo}
\end{figure}
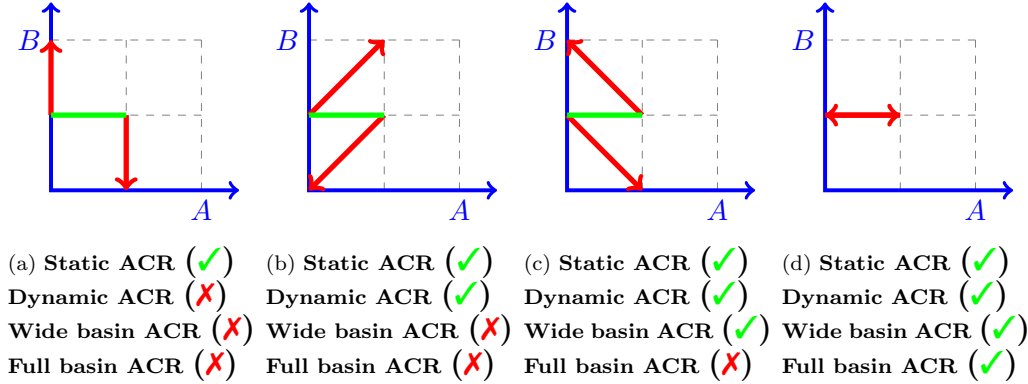

\begin{proposition} \label{prop:o3gh3o5hosi}
A network $\GG$ is a static (dynamic) ACR network if $\GG$ has a static (dynamic) ACR species. 
\end{proposition}
\begin{proof}
Suppose that $\GG$ has a static (dynamic) ACR species, label it by $X$. So, the concentration of $X$ is a static (dynamic) ACR variable for all choices of rate constants $K$. Then $(\GG,K)$ is a static (dynamic) ACR system for all choices of $K$, and so $\GG$ is a static (dynamic) ACR network. 
\end{proof}

The converse to Proposition \ref{prop:o3gh3o5hosi} may not hold. 
To see this, consider a static (dynamic) ACR network $\GG$. Let $K$ and $K'$ be some distinct choices of kinetics, so that $(\GG,K)$ and $(\GG,K')$ are both static (dynamic) ACR systems. Suppose that both systems have only one ACR variable, but these are concentrations of different species $X$ and $X'$ in cases of $(\GG,K)$ and $(\GG,K')$, respectively. Then $\GG$ does not have any ACR species. Whether such an example can be constructed within the mass action framework is an open question.

\begin{enumerate}[label={\bf Ex \arabic*}.,wide, labelwidth=!, labelindent=0pt]
\setcounter{enumi}{\value{break}}

\item\label{ex17} {\em (Dynamic ACR network with dynamic ACR species, has capacity for static ACR but no static ACR species)}
Consider the reaction network
\begin{align} \label{eq:p8hinohohu}
0 \stackrel[k_2]{k_1}{\rlas} A, \quad A+B  \xrightarrow{k_3} A, \quad B  \xrightarrow{k_4} 2B. 
\end{align}
The mass action ODE system is
\begin{align} \label{eq:245yoeirghg}
\dot a = k_2(k_1/k_2- a), \quad \dot b = -k_3 b(a - k_4/k_3). 
\end{align}
It is clear that $a \xrin k_1/k_2$ for any initial value in $\R^2_{\ge 0}$, and so $A$ is a dynamic ACR species with dynamic ACR value of $k_1/k_2$. There are no positive steady states if $k_1/k_2 \ne k_4/k_3$, and so for any such rate constants the system is not static ACR. This implies that there is no static ACR species. Note that for the special case $k_1/k_2 = k_4/k_3$, $a$ is a static ACR variable with static ACR value of $k_1/k_2$, which is the same as its dynamic ACR value. 
\setcounter{break}{\value{enumi}}
\end{enumerate}

A concrete representation of networks with various ACR properties is in Figure \ref{fig:24958heghuwo}, where we depict specific stoichoimetries by embedding the network in the plane. 
All networks in Figure \ref{fig:24958heghuwo} are static ACR.
By embedding the networks as close to the origin as possible, and making the reaction arrows as small as possible, we get minimal motifs of various ACR types: static ACR, narrow basin dynamic ACR, wide basin dynamic ACR, and full basin dynamic ACR.  

Based only on the examples in Figure \ref{fig:24958heghuwo}, it might be tempting to think that the reason \ref{fig:24958heghuwo}(b) is narrow basin dynamic ACR is because it is not mass conserving, unlike \ref{fig:24958heghuwo}(c) and \ref{fig:24958heghuwo}(d). However, there do exist mass conserving narrow basin dynamic ACR systems as the following example shows.

\begin{enumerate}[label={\bf Ex \arabic*}.,wide, labelwidth=!, labelindent=0pt]
\setcounter{enumi}{\value{break}}

\item\label{ex19} {\em (Existence of mass conserving, narrow basin dynamic ACR system/network)}
Consider the following reaction network
\begin{align} \label{eq:35pyhrhhwe}
X + 2Z \xrightarrow{k_1} 2X + Y, \quad X + Y + 2Z  \xrightarrow{k_2} 4Z. 
\end{align}
The mass action ODE system is
\begin{align} \label{eq:p25yidaippjpo}
\dot x = xz^2(k_1 - k_2 y), \quad \dot y = xz^2(k_1 - k_2 y), \quad \dot z = -2 xz^2(k_1 - k_2 y). 
\end{align}
A mass conservation law involving all the species is $x + y + z = c >0$, i.e. the sum of concentrations of all the species is constant over time. 
Clearly $Y$ is a static and a dynamic ACR species, i.e. for all rate constants, $y$ is a static and dynamic ACR variable with ACR value $k_1/k_2$. This means that, if an initial condition is compatible with the set $\{y^* = k_1/k_2\}$, then $y \xrin k_1/k_2$. 
\begin{claim}
For the mass action system \eqref{eq:p25yidaippjpo}, $\{(x,y,z)~|~y-x>k_1/k_2\}$ is not compatible with $\{(x,y,z)~|~y = k_1/k_2\}$. 
\end{claim}
\begin{proof}
Two points $(x_1,y_1,z_1)$ and $(x_2,y_2,z_2)$ are compatible with each other if and only if $(x_1-x_2,y_1-y_2,z_1-z_2) \in \spn(1,1,-2)$. 
Suppose that $(x_1,y_1,z_1) \in \{(x,y,z)~|~y-x>k_1/k_2\}$ and let $c_1  \coloneqq  y_1 - x_1 > k_1/k_2$, so that $(x_1,y_1,z_1) = (x_1,x_1 + c_1,z_1)$. Suppose that $(x_2,y_2,z_2) \in \{(x,y,z)~|~y  = k_1/k_2\}$, so that $(x_2,y_2,z_2) = (x_2,k_1/k_2,z_2)$. Then, 
\[
(x_1 - x_2 ,y_1 - y_2,z_1 - z_2) = (x_1 - x_2 , x_1 + c_1 - k_1/k_2, z_1 - z_2). 
\]
But then, $(y_1 - y_2) - (x_1 - x_2) = c_1 - k_1/k_2  + x_2 > x_2 > 0$, i.e. $(y_1 - y_2) \ne (x_1 - x_2)$, and so $(x_1-x_2,y_1-y_2,z_1-z_2) \notin \spn(1,1,-2)$. This proves the claim. 
\end{proof}
In particular, the set $C\coloneqq \{(0,k_1/k_2 + \alpha, 0) : \alpha \in \R_{> 0}\}$ is not compatible with $\{(x,y,z)~|~y = k_1/k_2\}$ and clearly $\{y: (x,y,z) \in C\}$ has no upper bound. 
This implies that $Y$ is a narrow basin dynamic ACR species in the network \eqref{eq:35pyhrhhwe}. 
\setcounter{break}{\value{enumi}}
\end{enumerate}

\section{ACR in Complex Balanced Systems} \label{sec:5yiu6fikfyuerlkgfye}

Complex balanced systems (as well as the more restrictive detailed balanced systems \cite{horn1972general}) have played a central role in study of mathematical models of reaction networks \cite{yu2018mathematical}. It turns out that a complete characterization of static and dynamic ACR property can be  found for complex balanced systems, as described below. 
While large families of static ACR systems are easy to come by, here we describe a large family of dynamic ACR systems. 
We give network conditions that guarantee dynamic ACR as well as network conditions that forbid dynamic ACR in complex balanced systems. 
\begin{theorem} \label{thm:3o8y02834th}
Suppose that $\GG$ is a reaction network such that for the choice of mass action kinetics $K$, the system $(\GG,K)$ is complex balanced. 
Let $\SS$ be the stoichiometric subspace of $\GG$. Let $e_i$ denote the standard basis unit vector, with $1$ in the $i$th component and $0$ elsewhere. 
The following statements are equivalent. 
\been[label={A}{\arabic*}.] 
\item $x_i$ is a static ACR variable. 
\item $e_i \in \SS$. 
\item There exist reactions  $y_1 \to y'_1,  y_2 \to y'_2, ..., y_m \to y'_m$ in $\GG$ such that for some $\lambda_i \in \Z$ we have $\sum_{i=1}^m \lambda_i (y'_i - y_i) =  \alpha X_i$, where $\alpha \in \Z \setminus \{0\}$. 
\enen
The following statements are equivalent. 
\been[label={B}{\arabic*}.]
\item $x_i$ is a dynamic ACR variable. 
\item $x_i$ is a wide basin dynamic ACR variable. 
\item $x_i$ is a full basin dynamic ACR variable. 
\enen
Moreover, Bi $\implies$ Aj  for $i,j \in \{1,2,3\}$. 
Furthermore, if the steady states of $(\GG,K)$ are globally attracting then all six statements are equivalent. 
\end{theorem}
\begin{proof} The equivalence of conditions {\em A2 } and {\em A3 } follows just from the definition of the 
stoichiometric subspace $\SS$.  
We now show that conditions {\em A1} and {\em A2} are also equivalent. 
Let $\wt x$ be a complex balanced steady state of $(\GG,K)$. 
Then from \cite{horn1972general}, the set $Z$ of all positive steady states of $(\GG,K)$ satisfies:
\[
\log Z = \log \wt x + \SS^\perp, 
\]
where $\SS^\perp$ is the orthogonal complement of $\SS$. 
Furthermore, the condition {\em A1} is equivalent to 
$
Z \subset \{ x_i = a_i^*\}
$
for some $a_i^* > 0$, which in turn is equivalent to $\log Z \subset \{y: y_i = \log a_i^* \}$, i.e. $\log Z$ is in a hyperplane parallel to a coordinate hyperplane. 
Since $\SS^\perp$ is a subspace and a translation of $\log Z$,  
\[
\SS^\perp \subset \{ x_i = 0 \},
\]
which in turn is equivalent to condition {\em A2}.

With regard conditions {\em B1, B2, B3}, it is clear that {\em B3} $\implies$ {\em B2} $\implies$ {\em B1}. To see that {\em B1} $\implies$ {\em B3}, note that {\em B1} $\implies$ {\em A1} $\implies$ {\em A2}. But, this implies that any initial value $x_0 \in \mathbb R^n_{>0}$ is compatible with any hyperplane of the form  $\{x \in \R^n_{>0} ~|~ x_i = a_i^*\}$ with $a_i^*>0$. Therefore, if $x_i$ is a dynamic ACR variable, then the convergence of solutions to the hyperplane $\{x \in \R^n_{>0} ~|~ x_i = a_i^*\}$  holds for any positive initial value, which implies {\em B3}.
\end{proof}
\begin{corollary}
Suppose that $\GG$ is a reaction network such that for the choice of mass action kinetics $K$, the system $(\GG,K)$ is complex balanced. If two complexes in the same linkage class differ only in the  species $X_i$,  then $x_i$ is a static ACR variable. In particular, if  $0 \to X_i$ or $X_i \to 0$ is in $\GG$, then  $x_i$ is a static ACR variable.
\end{corollary}

Since any weakly reversible network with deficiency $\delta = 0$ is complex balanced for any choice of rate constants, we also obtain the following.

\begin{corollary} \label{cor:3otthohwfl}
Suppose that $\GG$ is weakly reversible and has $0$ deficiency. 
Let $\SS$ be the stoichiometric subspace of $\GG$. Let $e_i$ denote the standard basis unit vector, with $1$ in the $i$th component and $0$ elsewhere. 
The following are equivalent:
\been[label={A}{\arabic*}.]
\item $X_i$ is a static ACR species. 
\item $e_i \in \SS$. 
\item There exist reactions  $y_1 \to y'_1,  y_2 \to y'_2, ..., y_m \to y'_m$ in $\GG$ such that for some $\lambda_i \in \Z$ we have $\sum_{i=1}^m \lambda_i (y'_i - y_i) =  \alpha X_i$, where $\alpha \in \Z \setminus \{0\}$. 
\enen
Moreover, the following properties are also equivalent:
\been[label={B}{\arabic*}.]
\item $X_i$ is a dynamic ACR species. 
\item $X_i$ is a wide basin dynamic ACR species. 
\item $X_i$ is a full basin dynamic ACR species. 
\enen
Moreover, Bi $\implies$ Aj  for $i,j \in \{1,2,3\}$. 
Furthermore, if the steady states of $(\GG,K)$ are globally attracting for any $K$, then all six properties above are equivalent. 
\end{corollary}
\begin{remark}[Regarding equivalence of {\em A} and {\em B} statements] \label{rem:q3ihoiqhgo}
It is known that the positive steady states of a complex balanced system $(\GG,K)$ are globally attracting (within the set of positive compatible points) for any $K$ if $\GG$ satisfies any of the following conditions.
\been
\item $\GG$ is connected \cite{anderson2011proof}. 
\item $\SS$ is at most three dimensional \cite{pantea2012persistence,craciun2013persistence}. 
\enen
According to the Global Attractor Conjecture \cite{craciun2009toric, craciun2015toric}, {\em any} complex balanced system is globally attracting. 
\end{remark}
\begin{corollary} \label{cor:qeorhorgh5697iu}
Suppose that a mass action system $(\GG,K)$ is mass conserving and complex balanced. Then $(\GG,K)$ is neither static ACR nor dynamic ACR. 
\end{corollary}
\begin{corollary} \label{cor:eogho3gh3oqgih}
If a network $\GG$ is deficiency zero, weakly reversible, and mass conserving, then $\GG$ does not have capacity for static ACR nor for dynamic ACR. 
\end{corollary}
\begin{enumerate}[label={\bf Ex \arabic*}.,wide, labelwidth=!, labelindent=0pt]
\setcounter{enumi}{\value{break}}
\item\label{ex203} {\em (Network that has no capacity for static ACR nor for dynamic ACR)}
Consider the reaction network $A \rlas B$. By Corollary \ref{cor:eogho3gh3oqgih}, the network does not have the capacity for static ACR nor for dynamic ACR. 
\setcounter{break}{\value{enumi}}
\end{enumerate}

For other conditions that thwart static ACR, see Theorem 9.7.1 in \cite{feinberg2019foundations}. 
\begin{enumerate}[label={\bf Ex \arabic*}.,wide, labelwidth=!, labelindent=0pt]
\setcounter{enumi}{\value{break}}
\item \label{ex1800} {\em (Network that is both static and dynamic ACR)} 
Consider the weakly reversible, deficiency zero reaction network $\GG$ depicted below.
\begin{align*}
  \begin{tikzpicture}[baseline={(current bounding box.center)}, scale=1]
   \node[state] (C)  at (0,0)  {$C$};
   \node[state] (AB)  at (2,1)  {$A+B$};
   \node[state] (B)  at (4,0)  {$B$};
   \path[->,scale=4]
    (C) edge[] node {} (AB)
    (AB) edge[] node {} (B)
    (B) edge[] node {} (C)
;
  \end{tikzpicture}
\end{align*}
By Corollary \ref{cor:3otthohwfl}, $A$ is a static ACR species. Furthermore, by Remark \ref{rem:q3ihoiqhgo} (1. or 2.), $A$ is a full basin dynamic ACR species.

Assuming the global attractor conjecture \cite{craciun2009toric, craciun2015toric} is true, we have the following fact: if a complex balanced system has full dimension, then it is dynamic ACR in all species. The following example illustrates this fact. 
\item \label{ex18} {\em (Network that is both static and dynamic ACR in all species)} 
Consider the weakly reversible, deficiency zero reaction network $\GG$ depicted below.
\begin{align*}
  \begin{tikzpicture}[baseline={(current bounding box.center)}, scale=1]
   \node[state] (A)  at (0,0)  {$A$};
   \node[state] (AB)  at (2,1)  {$A+B$};
   \node[state] (B)  at (4,0)  {$B$};
   \path[->,scale=4]
    (A) edge[] node {} (AB)
    (AB) edge[] node {} (B)
    (B) edge[] node {} (A)
;
  \end{tikzpicture}
\end{align*}
By Corollary \ref{cor:3otthohwfl}, both $A$ and $B$ are static ACR species. Furthermore, by Remark \ref{rem:q3ihoiqhgo} (1. or 2.), both $A$ and $B$ are full basin dynamic ACR species. 
\end{enumerate}

\section{Discussion and future work}

Biochemical reaction networks need to maintain robustness in their outputs against highly variable protein or enzyme concentrations. An example of this can be found in bacterial two-component signaling systems, a class that encompasses several thousands of systems \cite{alon2019introduction}. 
We refer to {\em empirical robustness} as the property that the long-term measured concentration of a biochemical species (say $X$) is independent of initial conditions of all species as well as independent of long-term values of other species besides $X$. 
In this paper, we have developed a mathematical framework which will allow proper modeling of empirical robustness. 
We refer to the mathematical property by dynamic ACR -- {\em there is a positive translation of a coordinate hyperplane that attracts all trajectories}.  
This single notion covers a wide variety of dynamical behaviors. For instance there could be globally attracting steady states on the ACR hyperplane or otherwise there could be attracting limit cycle oscillations confined entirely to the ACR hyperplane. 
Most intriguingly, the ACR hyperplane could be an attractor to unbounded trajectories, which means that unbounded trajectories nevertheless have a bounded and robust component. This has implications for robustness in growing systems.

The previous approach to model empirical robustness was to define static ACR, wherein all steady states are confined to the ACR hyperplane. The obvious problem with this approach is that static ACR by no means guarantees dynamic convergence to the steady states or to the ACR hyperplane. So the long-term behavior of static ACR systems may not show concentration robustness. 
Dynamic ACR not only remedies this problem, but has other advantages. 
As discussed in Proposition \ref{prop:simpleunion}, unlike static ACR, the property of dynamic ACR is unaffected when taking a union of networks with unrelated species, thus networks with dynamic ACR are structurally robust. In Theorem \ref{thm:q035yhtihgje}, we discuss connections between the dynamic ACR property and globally attracting steady states.

The range of dynamical behaviors captured by dynamic ACR suggests that finding network conditions for dynamic ACR will require a complex research program. In this paper, we gave necessary and sufficient conditions for dynamic ACR in complex balanced systems, an importance class of systems in reaction network theory. 
In \cite{joshi2022motifs}, we study static and dynamic ACR in small reaction networks, specifically those with 2 reactions and at most 2 species. 
We show that for such reaction networks, the network motif (the relative locations of the reactant complexes and the orientation of the reaction arrows in the Euclidean plane) is deeply connected with the dynamical properties, especially those related to ACR. 

Finally, in this paper, we have discussed consequences of static and dynamic ACR. Surprisingly, both static  and dynamic ACR are compatible not only with steady states but also with oscillations (\ref{ex10} and \ref{ex12}). In future work, we discuss biochemically realistic networks and some remarkable robustness properties of networks with dynamic ACR.

\subsection*{Acknowledgments}
BJ was supported by an IINA award from the CSUSM Advisory Council.
GC acknowledges support from NSF grant DMS-1816238 and from a Simons Foundation Fellowship. 
We thank the referees for careful reading and helpful comments.

\bibliographystyle{unsrt}
\bibliography{acr}

\end{document}